\documentclass[11pt]{article}
\usepackage{mathtools,amsmath,amsthm,amssymb,cases}
\usepackage{epsfig}
\usepackage{leftidx}
\usepackage{graphicx}
\usepackage{amsfonts}
\usepackage{color, epstopdf}
\usepackage{cite}
\usepackage{graphicx,wrapfig}
\usepackage{flafter}
\usepackage{fancyhdr}
\usepackage{stmaryrd}
\usepackage{graphicx,subfigure}
\usepackage{multicol}
\usepackage{dsfont}
\usepackage{booktabs,threeparttable}
\usepackage[center]{caption2}
\usepackage{epstopdf}
\usepackage [latin1]{inputenc}
\usepackage{multirow}
\usepackage{enumerate}
\usepackage{enumitem}
\usepackage{algpseudocode,algorithm,algorithmicx}

\newtheorem{theorem}{Theorem}[section]
\newtheorem{lemma}{Lemma}[section]
\newtheorem{proposition}{Proposition}[section]
\newtheorem{corollary}{Corollary}[section]
\newtheorem{example}{Example}
\newtheorem{remark}{Remark}

\newcommand{\zd}{\,\mathrm{d}}

\newcommand{\abs}[1]{\left|#1\right|}

\newcommand{\bra}[1]{\left(#1\right)}
\newcommand{\brab}[1]{\big(#1\big)}
\newcommand{\braB}[1]{\Big(#1\Big)}
\newcommand{\brat}[1]{(#1)}
\newcommand{\kbra}[1]{\left[#1\right]}
\newcommand{\kbrab}[1]{\big[#1\big]}
\newcommand{\kbraB}[1]{\Big[#1\Big]}
\newcommand{\kbrat}[1]{[#1]}
\newcommand{\myinner}[1]{\left\langle#1\right\rangle}
\newcommand{\myinnert}[1]{\langle#1\rangle}
\newcommand{\myinnerb}[1]{\big\langle#1\big\rangle}
\newcommand{\myinnerB}[1]{\Big\langle#1\Big\rangle}
\newcommand{\mynorm}[1]{\left\|#1\right\|}
\newcommand{\mynormb}[1]{\big\|#1\big\|}

\def\lan#1{\textcolor{blue}{#1}}

\title{Average energy dissipation rates of explicit exponential 
	Runge-Kutta methods for gradient flow problems}
\author{Hong-lin Liao
	\thanks{ORCID 0000-0003-0777-6832. School of Mathematics,
		Nanjing University of Aeronautics and Astronautics,
		Nanjing 211106, China; Key Laboratory of Mathematical Modeling
		and High Performance Computing of Air Vehicles (NUAA), MIIT, Nanjing 211106, China. Emails: liaohl@nuaa.edu.cn and liaohl@csrc.ac.cn.
		This author's work is supported by NSF of China	under grant number 12071216.}
	    \quad Xuping Wang\thanks{School of Mathematics, Nanjing University of Aeronautics and Astronautics, Nanjing 211106, China. Email: wangxp@nuaa.edu.cn}
}
\date{}

\begin{document}
	
	\maketitle
	
	\begin{abstract}
		We propose a unified theoretical framework to examine the energy dissipation properties at all stages of explicit exponential Runge-Kutta (EERK) methods for gradient flow problems. The main part of the novel framework is to construct the differential form of EERK method by using the difference coefficients of method and the so-called discrete orthogonal convolution kernels. As the main result, we prove that an EERK method can preserve the original energy dissipation law unconditionally if the associated differentiation  matrix is positive semi-definite. A simple indicator, namely average dissipation rate, is also introduced for these multi-stage methods to evaluate the overall energy dissipation rate of an EERK method such that one can choose proper parameters in some parameterized EERK methods or compare different kinds of EERK methods. Some existing EERK methods in the literature are evaluated from the perspective of preserving the original energy dissipation law and the energy dissipation rate.
    \lan{Some numerical examples are also included to support our theory.}
		\\[1ex]		
			\\[1ex]
	\end{abstract}

\section{Introduction}
\setcounter{equation}{0}
We propose a unified theoretical framework to examine the energy dissipation properties at all stages of  explicit exponential Runge-Kutta (EERK) methods for solving the following semi-discrete semilinear parabolic problem 
\begin{align}\label{problem: autonomous}
	u_h'(t)+L_hu_h(t)=g_h(u_h(t)),\quad u_h(t_0)=u_h^0,
\end{align}
where $L_h$ is a symmetric, positive definite matrix resulting from certain spatial discretization of stiff term, typically the Laplacian operator $-\Delta$ with periodic boundary conditions, and $g_h$ represents a nonlinear but non-stiff term. 
Without losing the generality, the finite difference method is assumed to approximate spatial operators and we define the discrete $L^2$ inner product $\myinner{u,v}:=v^Tu$ and the $L^2$ norm $\mynorm{v}:=\sqrt{\myinner{v,v}}$. Assume that there exists a non-negative Lyapunov function $G_h$ such that $g_h(v)=-\frac{\delta}{\delta v} G_h(v).$
Then the problem \eqref{problem: autonomous} can be formulated into a gradient system
\begin{align}\label{problem: gradient flows}
	\frac{\zd u_h}{\zd t}=-\frac{\delta E}{\delta u_h}\quad\text{with}\quad
	E[v_h]:=\frac1{2}\myinnert{v_h,L_hv_h}+\myinnert{G_h(v_h),1}.
\end{align}
The dynamics approaching the steady-state solution $u_h^*$, that is $L_hu_h^*=g_h(u_h^*)$, of this dissipative system \eqref{problem: autonomous} satisfies the following \textit{original} energy dissipation law
\begin{align}\label{problem: energy dissipation law}
	\frac{\zd E}{\zd t}=\myinnerB{\frac{\delta E}{\delta u_h},\frac{\zd u_h}{\zd t}}=
	-\myinnerB{\frac{\zd u_h}{\zd t},\frac{\zd u_h}{\zd t}}\le0.
\end{align}


In simulating the semilinear parabolic problems \eqref{problem: autonomous} 
and related gradient flow problems \eqref{problem: gradient flows}, 
explicit exponential (including exponential integrating factor and exponential time differencing) integrators turned out to be very competitive, see  \cite{CanoMoreta:2024,CelledoniMarthinsenOwren:2003,CoxMatthews:2002JCP,DimarcoPareschi:2011,DuJuLiQiao:2019SINUM,DuJuLiQiao:2021SIREV,HochbruckOstermann:2010ActaNu,JuLiQiaoZhang:2018MCOM,
	JuZhangZhuDu:2015,JuLiQiao:2022,Krogstad:2005JCP,KassamTrefethen:2005, WangJuDu:2016}. 
For a detailed overview of such integrators and
their implementation, we refer to \cite{DuJuLiQiao:2021SIREV,HochbruckOstermann:2010ActaNu,JuZhangZhuDu:2015,FasiGaudreaultLundSchweitzer:2024}. The main idea behind these methods is to treat the linear part of problem exactly and the nonlinearity in an
explicit way and dates back to the 1960s, see  \cite{CelledoniMarthinsenOwren:2003,CoxMatthews:2002JCP,Lawson:1967,Pope:1963,
	StrehmelWeiner:1992book,HochbruckLubichSelhofer:1998,Verwer:1977}. 
For stiff problems, Hochbruck and Ostermann \cite{HochbruckOstermann:2005SINUM} constructed explicit exponential Runge-Kutta (also called exponential time differencing Runge-Kutta, ETDRK)  methods with stiff orders up to four and established the convergence in an abstract Banach
space framework of sectorial operators and locally Lipschitz continuous nonlinearities. Luan and Ostermann \cite{LuanOstermann:2014} showed that there does not exist an EERK method of order five with less than or equal to six stages and constructed a fifth-order method with eight stages for semilinear parabolic
problems. For the stability properties of EERK
methods, Maset and Zennaro \cite{MasetZennaro:2009MCOM} derived sufficient conditions of unconditional contractivity and unconditional asymptotic stability and investigated some popular EERK methods with respect to the two stability properties.

In the past decade, the explicit ETDRK methods \cite{CelledoniMarthinsenOwren:2003,CoxMatthews:2002JCP,HochbruckOstermann:2010ActaNu}
became popular in simulating gradient flow problems, see  \cite{DuJuLiQiao:2019SINUM,DuJuLiQiao:2021SIREV,JuLiQiaoZhang:2018MCOM,
	JuZhangZhuDu:2015,JuLiQiao:2022,LiLiJuFeng:2021SISC,LiuQuanWang:2023,WangJuDu:2016,ZhangLiuQianSong:2023},
	in the context of partial differential equations. 
 One of main concerns is whether these ETD type methods can preserve the decaying of original energy $E[u_h(t)]$. Although the first-order ETD1 method has been proven in \cite{DuJuLiQiao:2019SINUM,DuJuLiQiao:2021SIREV,JuLiQiaoZhang:2018MCOM} 
 to preserve the energy decaying, the energy dissipation property of high-order EERK methods seems theoretically challenging due to their multi-stage nature.
 Very recently, the second-order ETD2RK method in \cite{CoxMatthews:2002JCP} has been shown to preserve the original energy decaying of the scalar gradient system \cite{FuYang:2022JCP} and the matrix gradient system \cite{LiuQuanWang:2023}.  These works are theoretically interesting, while their analysis may be limited since the proofs for the energy decaying heavily rely on technical skills and would be difficult to extend for other ETD methods or general situations, such as the parameterized EERK methods constructed by Hochbruck and Ostermann \cite{HochbruckOstermann:2005SINUM,HochbruckOstermann:2010ActaNu}.

In this article, we will focus on whether and to what extent the multi-stage EERK methods preserve the original energy dissipation law \eqref{problem: energy dissipation law}. In the next section, a unified theoretical framework for the stage energy dissipation property of EERK methods is established by constructing the differential forms of EERK methods and a new concept, namely average dissipation rate, is introduced for these multi-stage methods to evaluate the overall energy dissipation rate of an EERK method such that one can choose proper parameters in some parameterized EERK methods or compare different kinds of EERK methods. Our main results are stated in Theorem \ref{thm: energy stability} and Lemma \ref{lemma: average dissipation rate}.

As applications of our theory, three parameterized second-order EERK methods, including the widespread ETD2RK scheme \cite{CoxMatthews:2002JCP} and the three-stage method by Strehmel and Weiner \cite{StrehmelWeiner:1992book},  are discussed in Section 3. Some popular methods are evaluated and suggested for practical numerical simulations, see Table \ref{table: EERK2 Energy dissipation property}, in which the abscissa choices  for the contractivity and the energy stability of three second-order EERK methods are summarized. Section 4 addresses four third-order EERK methods, including the ETD3RK \cite{CoxMatthews:2002JCP}, ETD2CF3 \cite{CelledoniMarthinsenOwren:2003} and two parameterized methods developed by Hochbruck and Ostermann \cite{HochbruckOstermann:2005SINUM}. Table \ref{table: EERK3 Energy dissipation property} collects some abscissa choices  for the energy stability of four third-order EERK methods. \lan{Numerical experiments are presented in Section 5 to support our theory.}
Short comments on four fourth-order EERK methods from \cite{CoxMatthews:2002JCP,HochbruckOstermann:2005SINUM,Krogstad:2005JCP,StrehmelWeiner:1992book} and some concluding remarks on the new theory are presented in the last section. 

\section{Stage energy laws of EERK methods}
\setcounter{equation}{0}

\subsection{General class of EERK methods}
Let $u_h^k$ be the numerical approximation of $u_h(t_k)$ at the grid point $t_k$ for $0\le k\le N$. To integrate the semilinear parabolic problem \eqref{problem: autonomous} from the discrete time $t_{n-1}$ ($n\ge1$) to the next grid point $t_n=t_{n-1}+\tau$, the construction of one-step EERK methods (typically, $\tau$ also represents a variable-step size) starts from the following variation-of-constants formula  
\begin{align*}
	u_h(t_{n-1}+\tau)=e^{-\tau L_h}u_h(t_{n-1})
	+\int_{0}^{\tau}e^{-(\tau-\sigma) L_h}g_h\kbra{u_h(t_{n-1}+\sigma)}\zd\sigma.
	\end{align*}
Let $U^{n,i}$ be the approximation of $u_h(t_{n-1}+c_{i}\tau)$ at the abscissas $c_1:=0$, $c_i\in(0,1]$ for $2\le i\le s$, and $c_{s+1}:=1$. By replacing $\tau$ by $c_i\tau$ to define the internal stages $t_{n-1}+c_{i}\tau$, one can construct the following general class of EERK methods:
\begin{subequations}\label{Scheme: general EERK}	
	\begin{align}	
	&~U^{n,1}=u_h^{n-1},	\\
	&~U^{n,i+1}=\chi_{i+1}(-{\tau}L_h)U^{n,1}
	+\tau\sum_{j=1}^{i}a_{i+1,j}(-{\tau}L_h)g_h(U^{n,j}),
	\quad\text{$1\le i\le s-1$,}\\
	&~U^{n,s+1}=\chi(-{\tau}L_h)U^{n,1}+\tau\sum_{j=1}^{s}b_{j}(-{\tau}L_h)g_h(U^{n,j}),\\
	&~u_h^{n}=U^{n,s+1}.	
	\end{align}
\end{subequations}
The method coefficients $\chi_i$, $\chi$, $a_{ij}$ and $b_{j}$ are constructed from linear combinations of the entire functions $\varphi_j(z)$ and scaled versions thereof. These functions are given by
\begin{align}\label{def: varphi}
	\varphi_0(z)=e^z\quad\text{and}\quad\varphi_j(z):=\int_0^1e^{(1-s)z}\frac{s^{j-1}}{(j-1)!}\zd s \quad\text{for $z\in\mathbb{C}$ and $j\geq 1,$}
\end{align}
which satisfy the recursion formula
\begin{align}\label{def: varphi procedure}	
	\varphi_{k+1}(z)=\frac{\varphi_k(z)-1/k!}{z} \quad \text{for $k\geq 0$.}	
\end{align}	
Here the involved matrix functions $\varphi_j(-{\tau}L_h)$ are defined on the spectrum of $-{\tau}L_h$, that is, the values  $\{\varphi_j(\lambda_k): 1 \leq k \leq \mathrm{dim}(-{\tau}L_h)\}$ exist, where $\lambda_k$ are the eigenvalues of $-{\tau}L_h$ and thus $\varphi_j(\lambda_k)$ are the eigenvalues of $\varphi_j(-{\tau}L_h)$. More properties on the matrix functions can be found in \cite[Theorem 1.13]{Higham:2008book}, and, typically in this article, $f(-{\tau}L_h)$ is a positive definite operator if the given entire function $f$ is positive.

 
Always we assume that $\chi_i(0)=1$ and $\chi(0)=1$ for consistency. This scheme \eqref{Scheme: general EERK} reduces 
to an explicit Runge-Kutta method with coefficients $a_{ij}:=a_{ij}(0)$ and $b_{j}:=b_{j}(0)$
if we put $L_h=0$. 
The latter method will be called the \textit{underlying explicit Runge-Kutta method}
 henceforth. 
We suppose throughout the paper that the underlying Runge-Kutta method satisfies
\begin{align*}
\sum_{j=1}^{s}b_{j}(0)=1\quad\text{and}\quad	
\sum_{j=1}^{i-1}a_{ij}(0)=c_i\quad\text{for $ i=1,2,\cdots, s$}, 
\end{align*}
which makes it invariant under the transformation of \eqref{Scheme: general EERK} to the non-autonomous system.
A desirable property of numerical methods is that they preserve equilibria $u_h^*$ of
\eqref{problem: gradient flows}. Requiring $U^{n,i}=u_h^n=u_h^*$ for all
$i$ and $n \geq 0$ immediately yields the necessary and sufficient conditions. It turns out
that the method coefficients have to satisfy
\begin{align}\label{cond: equilibria}
\sum_{j=1}^{s}b_{j}(z)=\frac{\chi(z)-1}{z}\quad\text{and}\quad	
	\sum_{j=1}^{i-1}a_{ij}(z)=\frac{\chi_i(z)-1}{z}\quad\text{for $i=1,2,\cdots, s.$}
\end{align}
Without further mention, we consider the methods with 
$\chi(z)=e^{z}$ and $\chi_i(z)=e^{c_iz}$ for $1\leq i\leq s$.
With the help of \eqref{cond: equilibria}, the functions $\chi_i$ and $\chi$ 
can be eliminated in \eqref{Scheme: general EERK}. The numerical scheme \eqref{Scheme: general EERK} then takes the form
\begin{subequations}\label{Scheme: general EERK2}	
	\begin{align}				&~U^{n,i+1}=U^{n,1}+\tau\sum_{j=1}^{i}a_{i+1,j}(-{\tau}L_h)\kbra{g_h(U^{n,j})-L_hU^{n,1}}	\quad\text{for $1\le i\le s-1$,}\\
		&~U^{n,s+1}=U^{n,1}+\tau\sum_{j=1}^{s}b_{j}(-{\tau}L_h)\kbra{g_h(U^{n,j})-L_hU^{n,1}}.
	\end{align}
\end{subequations}
To simplify our notations,  define
\begin{align}\label{cond: Matrix setting}
a_{s+1,j}(z):=b_j(z),\quad 1\leq j\leq s.
\end{align} 
Then the EERK method \eqref{Scheme: general EERK2} applying to \eqref{problem: autonomous} reads
	\begin{align}	\label{Scheme: general EERK3}				&~U^{n,i+1}=U^{n,1}+\tau\sum_{j=1}^{i}a_{i+1,j}(-{\tau}L_h)\kbra{g_h(U^{n,j})-L_hU^{n,1}}	\quad\text{for $1\le i\le s$.}
	\end{align}
Always, we assume that $a_{k+1,k}(z)\neq0$ for any $1\le k\le s$. 
The associated Butcher tableau reads, where we use the abbreviations $a_{ij}:=a_{ij}(-{\tau}L_h)$,
\begin{equation*}
	\begin{array}{c|ccccc}
		c_1 & 0 &  &  &  &   \\
		c_{2} & a_{21} & 0 &  &  &   \\
		c_{3} & a_{31} & a_{32} & 0 &  &   \\
		\vdots & \vdots & \vdots & \ddots & \ddots &  \\[2pt]
		c_{s} & a_{s,1} & a_{s,2} &  \cdots  & a_{s,s-1}   & 0 \\[2pt]
		\hline  & a_{s+1,1} & a_{s+1,2} & \cdots 
		& a_{s+1,s-1}  & a_{s+1,s}
\end{array}\quad.\end{equation*}

 \subsection{Our theoretical framework}

Motivated by Du et al. \cite{DuJuLiQiao:2019SINUM,DuJuLiQiao:2021SIREV}, 
we introduce the stabilized operators with a parameter $\kappa\ge0$,
\begin{align}\label{def: stabilized parameter}
	L_{\kappa}:=L_h+\kappa I\quad\text{and}\quad g_{\kappa}(u):=g_h(u)+\kappa u,
\end{align} 
such that the problem \eqref{problem: autonomous} becomes the stabilized version
\begin{align}\label{problem: stabilized version}
	u_h'(t)=-L_{\kappa}u_h(t)+g_{\kappa}(u_h),\quad u_h(t_0)=u_h^0.
\end{align}
Thus, applying \eqref{Scheme: general EERK3} to \eqref{problem: stabilized version}, 
we have the following EERK method
\begin{align}\label{Scheme: general EERK stabilized2}			
	U^{n,i+1}=U^{n,1}+\sum_{j=1}^{i}a_{i+1,j}(-{\tau}L_{\kappa})\kbra{{\tau}g_{\kappa}(U^{n,j})-{\tau}L_{\kappa}U^{n,1}}
	\quad\text{for $1\le i\le s$.}
\end{align}

To make our idea more concise, we assume further that the nonlinear function $g_h$ is Lipschitz continuous with a constant $\ell_{g}>0$, cf. \cite{StuartHumphries:1998} or the recent discussions in \cite[subsection 2.2]{FuYang:2022JCP}. 
In theoretical manner, the stabilization parameter $\kappa$ in \eqref{def: stabilized parameter} is chosen properly large $\kappa\ge 2\ell_g$, see Remark \ref{remark: stabilization parameter} for an alternative choice, to enhance the dissipation of linear part so that the nonlinear growth of $g_h$ can be formally controlled in the numerical analysis. In this sense,  if an EERK method is proven to maintain the original energy dissipation law \eqref{problem: energy dissipation law} unconditionally, we mean that this EERK method can be stabilized by setting a properly large parameter $\kappa$ (which may not be necessary in actual calculations). 
To derive the energy dissipation law of the general EERK method \eqref{Scheme: general EERK stabilized2}, we need the following result. The proof is standard and we include it for completeness. 

\begin{lemma}\label{lemma: origional energy derivation}
	If $g_h$ is Lipschitz-continuous with a constant $\ell_{g}>0$ and 
	$\kappa\ge2\ell_g$, then
	\begin{align*}
		\myinnerb{u-v,g_{\kappa}(v)-\frac12L_{\kappa} (u+v)} \leq E[v]-E[u],
	\end{align*}
	where the energy $E$ is defined in \eqref{problem: gradient flows}.
\end{lemma}
\begin{proof}Since $g_h$ is Lipschitz continuous, 
	\cite[Lemma 2.8.20]{StuartHumphries:1998} gives
	\begin{align}\label{ieq:  Lipschitz continuous}
		\myinnerb{u-v,g_h(v)}\le \myinnerb{G_h(v)-G_h(u), 1}+\ell_g\mynormb{u-v}^2.
	\end{align}
	It follows that
	\begin{align*}
		&\,\myinnerb{u-v,g_{\kappa}(v)-\frac{\kappa}2(u+v)}
		=\myinnerb{u-v,g_h(v)-\frac{\kappa}2 (u-v)}\\
		=&\,\myinnerb{u-v,g_h(v)}-\frac{\kappa}2\mynormb{u-v}^2
		\le \myinnerb{G_h(v)-G_h(u), 1}-\frac{1}2\brat{\kappa-2\ell_g}\mynormb{u-v}^2.
	\end{align*}
	Also, it is easy to know that
	\begin{align*}
		\myinnerb{u-v,\frac{\kappa}2(u+v)-\frac12L_{\kappa} (u+v)}=\frac12\myinnerb{v-u,L_h(u+v)}	=&\,\frac{1}2\myinnert{v,L_hv}-\frac{1}2\myinnert{u,L_hu}.
	\end{align*}
	Adding up the above two results yields the claimed inequality and completes the proof.
\end{proof}
 
    Our theoretical framework contains three main steps: 
    \begin{description}
    	\item[(1)]  \underline{Compute difference coefficients}: we introduce a class of difference coefficients, for $i=1,2,\cdots,s$,
    	\begin{align}\label{def: difference coefficients A}			
    		\underline{a}_{i+1,i}(z):=a_{i+1,i}(z)\quad\text{and}\quad
    		\underline{a}_{i+1,j}(z):=a_{i+1,j}(z)-a_{i,j}(z)\quad\text{for $1\le j\le i-1$.}
    	\end{align}
    	It is not difficult to derive from \eqref{Scheme: general EERK stabilized2} that
    	\begin{align}\label{Scheme: general EERK stabilized2II}			
    		\delta_{\tau}U^{n,i+1}=&\,\sum_{j=1}^{i}\underline{a}_{i+1,j}(-{\tau}L_{\kappa})\kbra{{\tau}g_{\kappa}(U^{n,j})-{\tau}L_{\kappa}U^{n,1}}
    		\quad\text{for $1\le i\le s$,}
    	\end{align}
    	where the (stage) time difference $\delta_{\tau}U^{n,i+1}:=U^{n,i+1}-U^{n,i}$ for $1\le i\le s$. The associated Butcher difference (Butcher-Diff) tableau reads 
    	\begin{equation*}
    		\text{Butcher-Diff tableau:}\quad\begin{array}{c|ccccc}
    			c_1 & 0 &  &  &  &   \\
    			c_{2} & \underline{a}_{21} & 0 &  &  &   \\
    			c_{3} & \underline{a}_{31} & \underline{a}_{32} & 0 &  &   \\
    			\vdots & \vdots & \vdots & \ddots & \ddots &   \\[2pt]
    			c_{s} & \underline{a}_{s,1} & \underline{a}_{s,2} &  \cdots  & \underline{a}_{s,s-1}   & 0 \\[2pt]
    			\hline  & \underline{a}_{s+1,1} & \underline{a}_{s+1,2} & \cdots 
    			& \underline{a}_{s+1,s-1}  & \underline{a}_{s+1,s}
    	\end{array}\quad.\end{equation*}
    \item[(2)] \underline{Determine DOC kernels and differential form}: we introduce the so-called discrete orthogonal convolution (DOC) kernels 
    $\underline{\theta}_{k,j}(z)$ with respect to the coefficient $\underline{a}_{ij}$, cf. \cite{LiLiao:2022,LiaoTangZhou:2024,LiaoZhang:2021}, 
    \begin{align}\label{eq: orthogonal procedureII}
    	\underline{\theta}_{k,k}(z):=\frac1{\underline{a}_{k+1,k}(z)}\quad\text{and}\quad
    	\underline{\theta}_{k,j}(z):=-\sum_{\ell=j+1}^{k}\underline{\theta}_{k,\ell}(z)
    	\frac{\underline{a}_{\ell+1,j}(z)}{\underline{a}_{j+1,j}(z)}
    	\quad\text{for $1\leq j\le k-1$.}
    \end{align}
    It is easy to check the following discrete orthogonal identity,
    \begin{align}\label{eq: orthogonal identityII}
    	\sum_{\ell=j}^{m}\underline{\theta}_{m,\ell}(z)\underline{a}_{\ell+1,j}(z)\equiv\delta_{mj}\quad \text{for $1\leq j\leq m\leq s$},
    \end{align}
    where $\delta_{mj}$ is the Kronecker delta symbol with $\delta_{mj}=0$ if $j\neq m$.
    Multiplying the above equation \eqref{Scheme: general EERK stabilized2II}	 by the DOC kernels (matrices) $\underline{\theta}_{k,i}(-{\tau}L_{\kappa})$,  and summing $i$ from 1 to $k$, one can apply the discrete orthogonal identity \eqref{eq: orthogonal identityII} to find that
    \begin{align*}	
    	\sum_{i=1}^{k}\underline{\theta}_{k,i}(-{\tau}L_{\kappa})\delta_{\tau}U^{n,i+1}
    	=&\,\sum_{i=1}^{k}\underline{\theta}_{k,i}(-{\tau}L_{\kappa})\sum_{j=1}^{i}\underline{a}_{i+1,j}(-{\tau}L_{\kappa})
    	\kbra{{\tau}g_{\kappa}(U^{n,j})-{\tau}L_{\kappa}U^{n,1}}\\    
    	=&\,\sum\limits_{j=1}^{k}\sum_{i=j}^{k}\underline{\theta}_{k,i}(-{\tau}L_{\kappa})\underline{a}_{i+1,j}(-{\tau}L_{\kappa})
    	\kbra{{\tau}g_{\kappa}(U^{n,j})-{\tau}L_{\kappa}U^{n,1}}\\
    	=&\,{\tau}g_{\kappa}(U^{n,k})-{\tau}L_{\kappa}U^{n,1}\\
    	=&\,{\tau}g_{\kappa}(U^{n,k})-{\tau}L_{\kappa}U^{n,k+1}
    	+{\tau}L_{\kappa}\sum_{\ell=1}^{k}\delta_{\tau}U^{n,\ell+1}
    	\quad\text{for $1\le k\le s$.}
    \end{align*}
    Thus we have an equivalent form (differential form) of 
    the EERK method \eqref{Scheme: general EERK stabilized2}
    \begin{align}\label{Scheme: DOC EERK stabilized2II}			
    	\sum_{\ell=1}^{k}d_{k\ell}(-{\tau}L_{\kappa})\delta_{\tau}U^{n,\ell+1}
    	=	{\tau}g_{\kappa}(U^{n,k})-\frac{\tau}2L_{\kappa}(U^{n,k+1}+U^{n,k})
    	\quad\text{for $1\le k\le s$,}
    \end{align}
    where the functions $ d_{k\ell}$ are defined by
    \begin{align}\label{Def: Differential Matrix DII}			
    d_{k\ell}(z):=\underline{\theta}_{k\ell}(z)+\frac{z}{2}\brab{2-\delta_{k\ell}}\quad\text{for $1\le \ell\le k\le s$}\quad\text{and}\quad 
    	d_{k\ell}(z):=0\quad\text{for $\ell> k$}.
    \end{align}
    The associated lower triangular matrix $D:=(d_{k\ell})_{s\times s}$  is called the differentiation  matrix.
Always, we denote the symmetric part $\mathcal{S}(D;z):=\frac{1}{2}\kbrat{D(z)+D(z)^T}.$ 
    
    \item[(3)] \underline{Establish stage energy dissipation law}: this process is standard and we have the following result,
    which simulates the original energy dissipation law \eqref{problem: energy dissipation law} at all stages.
    
    \begin{theorem}\label{thm: energy stability} If $\mathcal{S}(D;z)$ is positive (semi-)definite, 
    	the EERK method \eqref{Scheme: general EERK stabilized2} preserves
    	the original energy dissipation law \eqref{problem: energy dissipation law} 
    	at all stages without any time-step constraints, 
    	\begin{align}\label{thmResult: stage energy laws}			
    		E[U^{n,j+1}]-E[U^{n,1}]\le&\,-\frac1{\tau}\sum_{k=1}^{j}\myinnerB{\delta_{\tau}U^{n,k+1},
    			\sum_{\ell=1}^{k}d_{k\ell}(-{\tau}L_{\kappa})\delta_{\tau}U^{n,\ell+1}}
    			\quad\text{for $1\le j\le s$,}
    	\end{align}
    	and in particular, by taking $j:=s$,
    	\begin{align*}			
    		E[u_h^{n}]-E[u_h^{n-1}]\le-\frac1{\tau}\sum_{k=1}^{s}\myinnerB{\delta_{\tau}U^{n,k+1},
    			\sum_{\ell=1}^{k}d_{k\ell}(-{\tau}L_{\kappa})\delta_{\tau}U^{n,\ell+1}}\quad\text{for $n\ge1$.}
    	\end{align*}
    \end{theorem}
    \begin{proof}Making the inner product of the equivalent form \eqref{Scheme: DOC EERK stabilized2II} 
    	with $\frac1{\tau}\delta_{\tau}U^{n,k+1}$ and 
    	summing $k$ from $k=1$ to $j$, one can find that
    	\begin{align*}			
    		\frac1{\tau}\sum_{k=1}^{j}
    		\myinnerB{\delta_{\tau}U^{n,k+1},
    			\sum_{\ell=1}^{k}d_{k\ell}(-{\tau}L_{\kappa})\delta_{\tau}U^{n,\ell+1}}
    		=\sum_{k=1}^{j}\myinnerB{\delta_{\tau}U^{n,k+1},
    			g_{\kappa}(U^{n,k})-\frac{1}2L_{\kappa}(U^{n,k+1}+U^{n,k})}
    	\end{align*}
    	for $1\le j\le s$. Lemma \ref{lemma: origional energy derivation} yields
    	the following energy dissipation law at each stage
    	\begin{align*}			
    		E[U^{n,j+1}]-E[U^{n,1}]
    		+&\,\frac1{\tau}\sum_{k=1}^{j}\myinnerB{\delta_{\tau}U^{n,k+1},
    			\sum_{\ell=1}^{k}d_{k\ell}(-{\tau}L_{\kappa})\delta_{\tau}U^{n,\ell+1}}\le0
    	\end{align*}    	 
    	for $1\le j\le s$. It completes the proof.    	
    \end{proof}

    For $1\le j\le s$, let $D_{j}:=D[1:j,1:j]$ be the $j$-th sequential sub-matrix of the matrix $D$. 
    We denote further that  $\delta_{\tau}\vec{U}_{n,j+1}:=(\delta_{\tau}U^{n,2},\delta_{\tau}U^{n,3},\cdots,\delta_{\tau}U^{n,j+1})^T.$ The above stage energy dissipation law \eqref{thmResult: stage energy laws}
    can be formulated as 
    \begin{align}\label{thmResult: stage energy laws Matrix}			
    	E[U^{n,j+1}]-E[U^{n,1}]\le-\frac1{\tau}\myinnerB{\delta_{\tau}\vec{U}_{n,j+1},
    		D_{j}(-{\tau}L_{\kappa})\delta_{\tau}\vec{U}_{n,j+1}}
    	\quad\text{for $1\le j\le s$.}
    \end{align}    
    After the completion of this article, we are informed that, by computing the original energy difference $E[u_h^{n}]-E[u_h^{n-1}]$ with a key inequality, Fu, Shen and Yang independently derived the same sufficient condition of Theorem 2.1, cf. \cite[Theorem 2.1]{FuShenYang:2024}, to ensure that the EERK method \eqref{Scheme: general EERK stabilized2} maintains the decreasing of original energy, that is, $E[u_h^{n}]\le E[u_h^{n-1}]$. In the following subsection, we will introduce a simple indicator for evaluating to what extent the multi-stage EERK method \eqref{Scheme: general EERK stabilized2} preserves the original energy dissipation law \eqref{problem: energy dissipation law}.
    \end{description}
    
    \subsection{Averaged dissipation rate}
    Theorem \ref{thm: energy stability} shows that the EERK method \eqref{Scheme: general EERK stabilized2} is unconditionally energy stable if the differentiation  matrix $D(z)$ is semi-positive definite, that is, all eigenvalues $\lambda_{i}(z)$ $(i=1,2,\cdots,s)$ of the symmetric part $\mathcal{S}(D;z)$
    are nonnegative. A necessary condition is that the average eigenvalue is nonnegative, $$\overline{\lambda}(z):=\frac1{s}\sum_{i=1}^s\lambda_{i}(z)
    =\frac1{s}\mathrm{tr}\brab{D(z)}\ge0.$$
    If $\lambda_{\min}\le\lambda_{i}(z)\le \lambda_{\max}$ $(i=1,2,\cdots,s)$ for any $z\le0$, one has
    \begin{align*}			
    	\lambda_{\min}\myinnerb{\vec{v},\vec{v}}
    	\le\myinnerb{\vec{v},D(-{\tau}L_{\kappa})\vec{v}}
    	\le \lambda_{\max}\myinnerb{\vec{v},\vec{v}}.
    \end{align*} 
    Then, according to \eqref{thmResult: stage energy laws Matrix}, the overall 
    energy dissipation rate of the energy $E[u_h^{n}]$ could be roughly estimated by the average  
    eigenvalue $\overline{\lambda}(z)$ of $\mathcal{S}(D;z)$. If $\overline{\lambda}(z)\ge0$, one could use the following \textit{average dissipation rate}
    \begin{align}\label{def: numerical rate}			
    	\mathcal{R}(z):=\frac1{s}\mathrm{tr}\brab{D(z)}	\quad\text{for $z\le0$,}
    \end{align} 
    to examine the energy dissipation behaviors among different methods,     see detailed arguments for second-order EERK methods in the next section.
    By using the definitions \eqref{Def: Differential Matrix DII} 
    and \eqref{eq: orthogonal procedureII} to compute the diagonal elements $d_{kk}(z)$ for $1\le k\le s$, it is not difficult to obtain the following result. 
    
    \begin{lemma}\label{lemma: average dissipation rate} If the EERK method \eqref{Scheme: general EERK stabilized2} preserves  the original energy dissipation law \eqref{problem: energy dissipation law} unconditionally, then the average dissipation rate  is nonnegative, that is,
    	\begin{align*}			
    		\mathcal{R}(z)=\frac{z}{2}+\frac{1}{s}\sum_{i=1}^s\frac{1}{a_{i+1,i}(z)}\ge0\quad\text{for $z\le0$}.
    	\end{align*}    
    \end{lemma}
    
     Typically, if $\mathcal{R}(z)>1$, the discrete energy $E[u_h^{n}]$ decays faster than the continuous counterpart $E[u_h(t_n)]$ and the dynamics approaching the steady-state solution appears a time ``ahead" effect. If $0<\mathcal{R}(z)<1$, the discrete energy $E[u_h^{n}]$ may decay slower and the dynamics appears a time ``delay" effect. In general, a time-stepping method is a ``good" candidate to preserve 
     the original energy dissipation law \eqref{problem: energy dissipation law} unconditionally if the average dissipation rate $\mathcal{R}(z)$ is nonnegative for $z\le0$ and is as close to 1 as possible within properly large range of $z$. Lemma \ref{lemma: average dissipation rate} provides us a simple criterion to evaluate the overall energy dissipation rate of an EERK method and then choose proper parameters in some parameterized EERK methods or compare different EERK methods.

 \begin{remark}\label{remark: stabilization parameter}
 The differential form \eqref{Scheme: DOC EERK stabilized2II} and the associated differentiation  matrix $D(z)$ of the EERK method \eqref{Scheme: general EERK stabilized2} would be ``optimal'' to evaluate the energy dissipation property although they are not unique. A direct choice is to retain only the pure implicit form of stiff term, that is, 
 	\begin{align}\label{Scheme: DOC EERK stabilized2II remark}	
 		\sum_{\ell=1}^{k}\tilde{d}_{k\ell}(-{\tau}L_{\kappa})\delta_{\tau}U^{n,\ell+1}
 		=	{\tau}g_{\kappa}(U^{n,k})-\tau L_{\kappa}U^{n,k+1}
 		\quad\text{for $1\le k\le s$,}
 	\end{align}
 	where the elements of differentiation  matrix $\widetilde{D}:=(\tilde{d}_{k\ell})_{s\times s}$ are defined by \begin{align}\label{Def: Differential Matrix DII remark}			
 		\tilde{d}_{k\ell}(z):=\underline{\theta}_{k\ell}(z)+z\quad\text{for $1\le \ell\le k\le s$}\quad\text{and}\quad 
 		\tilde{d}_{k\ell}(z):=0\quad\text{for $\ell> k$}.
 	\end{align}
 	If $\mathcal{S}(\widetilde{D};z)$ is positive (semi-)definite, one can follow the proof of Theorem \ref{thm: energy stability} to get
 		\begin{align}\label{thmResult: stage energy laws remark}			
 		E[U^{n,j+1}]-E[U^{n,1}]\le&\,-\frac1{\tau}\sum_{k=1}^{j}\myinnerB{\delta_{\tau}U^{n,k+1},
 			\sum_{\ell=1}^{k}\tilde{d}_{k\ell}(-{\tau}L_{\kappa})\delta_{\tau}U^{n,\ell+1}}\nonumber\\
 		&\,-\frac{1}{\tau}\sum_{k=1}^{j}\myinnerB{\delta_{\tau}U^{n,k+1},\frac{1}{2}\tau L_{\kappa}\delta_{\tau}U^{n,k+1}}
 		\quad\text{for $1\le j\le s$,}
 	\end{align}
  	in which the following result similar to Lemma \ref{lemma: origional energy derivation} has been used, 
  	\begin{align*}
  		\myinnert{u-v,g_{\kappa}(v)-L_{\kappa} u} \leq E[v]-E[u]-\frac{1}{2}\myinnert{u-v,L_{\kappa}(u-v)}\quad\text{for $\kappa\ge\ell_g$.}
  	\end{align*}
  	  It is easy to see that the positive (semi-)definiteness of $\mathcal{S}(\widetilde{D};z)$ is much severer than the condition of Theorem \ref{thm: energy stability} because the energy dissipation estimate \eqref{thmResult: stage energy laws remark} ignores the dissipation effect of the last term compared with \eqref{thmResult: stage energy laws}. Correspondingly, the overall dissipation rate will be also underestimated via the average dissipation rate $\widetilde{\mathcal{R}}(z)$, that is,
  	  \begin{align*}			
  	  	\widetilde{\mathcal{R}}(z):=\frac1{s}\mathrm{tr}\brab{\widetilde{D}(z)}=z+\frac{1}{s}\sum_{i=1}^s\frac{1}{a_{i+1,i}(z)}<{\mathcal{R}}(z)\quad\text{for $z<0$}.
  	  \end{align*} 
  	  In this situation, one may make misjudgment on the energy dissipation property of EERK methods. For example,  consider the one-parameter EERK2 method \eqref{scheme: EERK2 Butcher tableau} described below with
  	  \begin{align*}
  	  	\widetilde{\mathcal{R}}(c_2,z):=z+\frac1{2c_{2}\varphi_{1}(c_2z)}+\frac{c_{2}}{2\varphi_{2}(z)}.
  	  \end{align*}
  	  It is easy to know that 
  	  $\lim\limits_{z\rightarrow-\infty}\widetilde{\mathcal{R}}(c_2,z)=-\infty$ if $c_2\in(0,1)$, while $\lim\limits_{z\rightarrow-\infty}\widetilde{\mathcal{R}}(1,z)=\frac12.$ 
  	 This directly leads to incorrect conclusion that the EERK2 method with $c_2=1$ is 
  	 the only possible case to preserve 
  	  the energy dissipation law \eqref{problem: energy dissipation law}; In contrast, Corollary \ref{corollary: ETDRK2} says that the EERK2 method preserves
  	  the energy dissipation law \eqref{problem: energy dissipation law} unconditionally for $c_2\in[\tfrac1{2},1]$. In summary, the condition of Theorem \ref{thm: energy stability} is nearly ``optimal'' although we can not claim that the positive semi-definiteness of differentiation  matrix $D(z)$ is also necessary to the energy stability of the EERK method (the only loss of dissipation rate comes from the inequality \eqref{ieq:  Lipschitz continuous}  for controlling the nonlinear growth).
 \end{remark}

\subsection{Simple case: ETD1}

To end this section, we consider a simple case with $s=1$. 
The only reasonable choice is the exponential forward Euler \cite{HochbruckOstermann:2005SINUM}
or ETD1 \cite{CoxMatthews:2002JCP,DuJuLiQiao:2019SINUM} method with stiff order one. 
Applied to \eqref{problem: stabilized version}, it is
\begin{align}\label{scheme: EERK1 HO-from}	
	\delta_{\tau}u^{n,2}=\varphi_1(-{\tau}L_{\kappa})\kbra{{\tau}g_{\kappa}(u^{n,1})-{\tau} L_{\kappa}u^{n,1}}
\end{align}
or, recalling the recursive formula \eqref{def: varphi procedure},
\begin{align}\label{scheme: EERK1 MBP-from}	
	u^{n,2}=\varphi_0(-{\tau}L_{\kappa})u^{n,1}+{\tau}\varphi_1(-{\tau}L_{\kappa})g_{\kappa}(u^{n,1}).
\end{align}
The associated Butcher and Butcher-Diff tableaux are the same, that is,
\begin{align*}
	\text{ETD1 Butcher or Butcher-Diff:}\quad\begin{array}{c|c}
	0 & 0    \\
	\hline  & \varphi_1 
\end{array}\;.
\end{align*}
The definition \eqref{Def: Differential Matrix DII} gives 
$$D^{(1)}=(d_{11}^{(1)})\quad\text{with}\quad d_{11}^{(1)}(z)=\frac{z}2+\frac1{\varphi_1(z)}=\frac{z(1+e^{-z})}{2(1-e^{-z})}\ge1\quad\text{for $z\le0$}.$$
Here and hereafter, the superscript $(p)$ is always used to indicate the order of the method, that is to say, $D^{(p)}$ and $\mathcal{R}^{(p)}$ denote the associated differential matrix and the average dissipation rate, respectively, of a formal $p$-th order EERK method.
Obviously, Theorem \ref{thm: energy stability} yields
\begin{corollary}\label{corollary: ETD1} 
	The exponential forward Euler \eqref{scheme: EERK1 HO-from}	preserves
	the energy dissipation law \eqref{problem: energy dissipation law},
	\begin{align*}			
		E[u_h^{n}]-E[u_h^{n-1}]\le-\frac{1}{\tau}\myinnerb{\delta_{\tau}u^{n},
			D^{(1)}(-{\tau}L_{\kappa})\brab{\delta_{\tau}u^{n}}}
		\quad\text{for $n\ge1$.}
	\end{align*}
\end{corollary}

By the definition \eqref{def: numerical rate}, one has $\mathcal{R}^{(1)}(z):=d_{11}^{(1)}(z)$ such that $\mathcal{R}^{(1)}(z)\ge1$ for any $z<0$ and $\lim_{z\rightarrow-\infty}\mathcal{R}^{(1)}(z)=+\infty$. It means that the dissipation rate of the discrete energy $E[u_h^{n}]$ approaches the original rate as the step size $\tau\rightarrow0$; while the exponential forward Euler \eqref{scheme: EERK1 HO-from} always generates a time ``ahead" (compared with the continuous counterpart $E[u_h(t_{n})]$) for any time-step sizes. 

By the form \eqref{scheme: EERK1 MBP-from}, it is easy to find that the ETD1 method is unconditionally contractive, also see \cite{MasetZennaro:2009MCOM}. The contractivity of EERK methods is essential to preserve the maximum bound principle of semilinear parabolic problems, cf. \cite{DuJuLiQiao:2019SINUM,DuJuLiQiao:2021SIREV,LiLiJuFeng:2021SISC,LiuQuanWang:2023,ZhangLiuQianSong:2023} and references therein; while detailed discussions are out of our current scope in this article.  


\section{Discrete energy laws of second-order methods}
\setcounter{equation}{0}

Second-order methods require two internal stages, $s=2$, at least. Hochbruck and Ostermann \cite{HochbruckOstermann:2005SINUM} derived the following stiff order conditions (the stiff order describes the behavior of the local error independently of the norm of the matrix $L_{\kappa}$) with a parameter $c_2$ $(0<c_2\le1)$ 
\begin{subequations}\label{conds: stiff two order}
	\begin{align}
		&~a_{31}(-{\tau}L_{\kappa})+a_{32}(-{\tau}L_{\kappa})=\varphi_1(-{\tau}L_{\kappa}),\label{6.12a}\\
		&~a_{32}(-{\tau}L_{\kappa})c_2=\varphi_2(-{\tau}L_{\kappa}),\label{6.12b}\\
		&~a_{21}(-{\tau}L_{\kappa})=c_2\varphi_1(-c_2{\tau}L_{\kappa}). \label{6.12c}
	\end{align}
\end{subequations}
They lead to the following one-parameter family of second-order EERK (EERK2) 
method with the following  Butcher tableau
\begin{align}\label{scheme: EERK2 Butcher tableau}
	\begin{array}{c|cccccc}
	0 &  &   \\
	c_{2} & c_{2}\varphi_{1,2} &   \\
	\hline  & \varphi_1-\frac{1}{c_2}\varphi_2 & \frac{1}{c_2}\varphi_2 
\end{array},
\end{align}
where the notations $\varphi_{i,j}$ are defined by 
\begin{align}\label{def: varphi ij abbreviated}
	\varphi_{i,j}:=\varphi_{i,j}(-{\tau}L_{\kappa})=\varphi_{i}(-c_j{\tau}L_{\kappa}),\quad i\ge0,\;1\leq j\leq s+1.
\end{align}
Note that, these abbreviations will be also used in the Butcher tableaus below.
This EERK2 method \eqref{scheme: EERK2 Butcher tableau} fulfills all conditions in \eqref{conds: stiff two order} 
and thus is stiff order two.  If the abscissa $c_2=1$,  it reduces to the so-called ETD2RK \cite{CoxMatthews:2002JCP,DuJuLiQiao:2019SINUM,FuYang:2022JCP} with the following form 
\begin{subequations}\label{scheme: ETDRK2}
	\begin{align}			
		U^{n,2}=&\,\varphi_{0}(-{\tau}L_{\kappa})U^{n,1}
		+{\tau}\varphi_{1}(-{\tau}L_{\kappa})
		g_{\kappa}(U^{n,1}),\\
		U^{n,3}=&\,U^{n,2}
		+{\tau}\varphi_2 (-{\tau}L_{\kappa})
		\kbra{g_{\kappa}(U^{n,2})-g_{\kappa}(U^{n,1})}.
	\end{align}
\end{subequations}
This case is also the only scenario to ensure the unconditional contractivity of EERK2 method \eqref{scheme: EERK2 Butcher tableau}, cf. \cite{MasetZennaro:2009MCOM}, due to the fact $\varphi_1(z)\ge\varphi_2(z)$ for $z\le0$.

By weakening the condition (\ref{6.12b}) to $a_{32}(0)c_2=\varphi_2(0)=\frac{1}{2}$, 
one has a one-parameter weak variant (called EERK2-w in short)
\begin{align}\label{scheme: EERK2-weak Butcher tableau}
	\begin{array}{c|cccccc}
		0 &  &   \\
		c_{2} & c_{2}\varphi_{1,2} &   \\
		\hline  & (1-\frac{1}{2c_2})\varphi_1 & \frac{1}{2c_2}\varphi_1
\end{array}\quad.
\end{align}
Although the EERK2-w method \eqref{scheme: EERK2-weak Butcher tableau}  does not have stiff order two, it achieves stiff convergence order two \cite[Section 5.1]{HochbruckOstermann:2005SINUM} under certain requirements (boundedness) on the discrete operator $\tau L_{\kappa}$. 
In the following, we consider their stage energy dissipation properties
in simulating  \eqref{problem: stabilized version}. 

\subsection{EERK2 method}
To establish the stage energy laws, we present the Butcher-Diff tableau
\begin{align}\label{scheme: EERK2 Butcher-Diff}
	\text{EERK2 Butcher-Diff:}\qquad
	\begin{array}{c|cccccc}
		0 &  &   \\
		c_2 & c_{2}\varphi_{1,2} &   \\
		\hline  & \varphi_1-\frac{1}{c_2}\varphi_2-c_{2}\varphi_{1,2} & \frac{1}{c_2}\varphi_2 
	\end{array}\quad.
\end{align}
By the procedure \eqref{eq: orthogonal procedureII}, 
one can compute the associated DOC kernels
\begin{align*}
	\underline{\theta}_{11}(z)=&\,\frac{1}{c_2\varphi_{1}(c_2z)},\quad 
	\underline{\theta}_{22}(z)=\frac{c_2}{\varphi_{2}(z)}\quad\text{and}\quad
	\underline{\theta}_{21}(z)
	=\frac{c_{2}\varphi_{1}(c_2z)-\varphi_1(z)+\frac{1}{c_2}\varphi_2(z)}{\varphi_{2}(z)\varphi_{1}(c_2z)}.
\end{align*}
The definition \eqref{Def: Differential Matrix DII} gives 
the following differentiation matrix
\begin{align*}
	D^{(2)}(c_2,z):=\begin{pmatrix}
		\frac{1}{c_2\varphi_{1}(c_2z)}+\frac{z}2 & 0  \\[4pt]	
		\frac{c_{2}\varphi_{1}(c_2z)-\varphi_1(z)+\frac{1}{c_2}\varphi_2(z)}{\varphi_{2}(z)\varphi_{1}(c_2z)}
		+z & \frac{c_2}{\varphi_{2}(z)}+\frac{z}2\\
	\end{pmatrix}.
\end{align*}
Now we consider the matrix $\mathcal{S}(D^{(2)};c_2,z)$, the symmetric part of $D^{(2)}(c_2,z)$.
The first leading principal minor  reads 
\begin{align*}
	\mathrm{Det}\kbrab{\mathcal{S}(D_1^{(2)};c_2,z)}=&\,
	d_{11}^{(2)}(c_2,z)=\frac{z(e^{c_2z}+1)}{2(e^{c_2z}-1)}\ge \frac1{c_2}\quad\text{for $c_2\in(0,1]$ and $z\le0$.}
\end{align*}
The second leading principal minor (determinant) of $\mathcal{S}(D^{(2)};c_2,z)$ is given by
\begin{align*}
	\mathrm{Det}\kbrab{\mathcal{S}(D^{(2)};c_2,z)}
	=&\,\frac{ (e^{c_2 z}-1)^{-2}z^2}{4(z-e^z+1)^2}g_{21}(c_2,z)\quad\text{for $z<0$,}
\end{align*}
where the auxiliary function $g_{21}$ is defined by
\begin{align}\label{def: g21} 
g_{21}(c_2,z):=&2 c_2 z e^{(c_2+2) z}-c_2^2 z^2 e^{2 c_2 z}
-2 c_2 z e^{c_2 z+z} (1-(c_2-1) z)-e^{2 z} (c_2^2 z^2+1)\nonumber\\
&\,+2e^z (1-(c_2-1) z)+(z+1) ((2 c_2-1) z-1)\quad\text{for $z<0$.}
\end{align}
Now we develop a technique of comparison function to handle the function $g_{21}$.
\begin{proposition}\label{proposition: g21}
	The function $g_{21}$ in \eqref{def: g21} is positive for  $c_2\in[\frac12,1]$ and $z<0$.
	\end{proposition}
\begin{proof}The condition $c_2\in[\frac12,1]$ comes from the simple fact $\lim_{z\rightarrow-\infty}g_{21}(c_2,z)/z^2=2c_2-1\ge0$.
	To handle $g_{21}$, we consider a comparison function (by setting $c_2:=\tfrac{1}{2}$ in all exponents of $g_{21}$)
	\begin{align*}
		g_{21}^*(c_2,z):=&2 c_2 z e^{\frac{5z}2}-c_2^2 z^2 e^{z}
		-2 c_2 z e^{\frac{3z}2} (1-(c_2-1) z)-e^{2 z} (c_2^2 z^2+1)\nonumber\\
		&\,+2e^z (1-(c_2-1) z)+(z+1) ((2 c_2-1) z-1)\quad\text{for $z<0$.}
	\end{align*}
	It is easy to check that $g_{21}(c_2,z)\ge g_{21}^*(c_2,z)$ for  $c_2\in[\frac12,1]$ and $z<0$. Actually,
	\begin{align*}
	g_{21}(c_2,z)-g_{21}^*(c_2,z)=&\,c_2 e^z (e^{z/2}-e^{c_2 z}) 
	\kbrab{c_2z^2 (e^{-z/2}+ e^{c_2 z-z}-2)-2z (e^z-1-z)}\\
	\ge&\,	c_2 e^z (e^{z/2}-e^{c_2 z}) 
	\kbrab{c_2z^2 (e^{-z/2}-1)-2z  (e^z-z-1)}\ge0
	\end{align*}
due to the facts $e^{-z/2}-1>0$ and $e^z-z-1>0$ for $z<0$.
	
    Note that, $g_{21}^*$ is a concave, quadratic polynomial with respect to $c_2$, that is,
    \begin{align*}
    	g_{21}^*(c_2,z)=-e^z(e^{z/2}-1)^2 z^2c_2^2+2z(1-e^{\frac{3 z}{2}})(z+1-e^z)c_2-(z+1-e^z)^2
    	\quad\text{for $z<0$.}
    \end{align*}
	Moreover, one can check that (technical details are omitted here), 
	cf. Figure \ref{fig: g21, g22 comparfuns}(a), 
	\begin{align*}		
		g_{21}^*(1,z)>0\quad\text{and}\quad
		g_{21}^*(\tfrac{1}{2},z)>0
		\quad\text{for $z<0$.}
	\end{align*}
	 They imply that $g_{21}^*(c_2,z)>0$ and then $g_{21}(c_2,z)>0$ for $c_2\in[\frac12,1]$ and $z<0$.
\end{proof}

\begin{figure}[htb!]
	\centering
	\subfigure[ $g^*_{21}(c_2,z)$]
	{\includegraphics[width=2.15in]{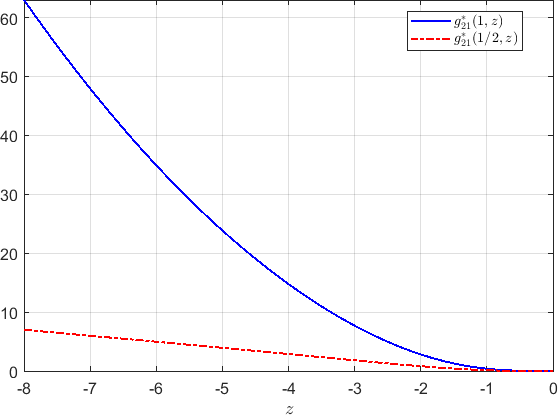}}\hspace{10mm}
	\subfigure[ $g^*_{22}(c_2,z)$]
	{\includegraphics[width=2.15in]{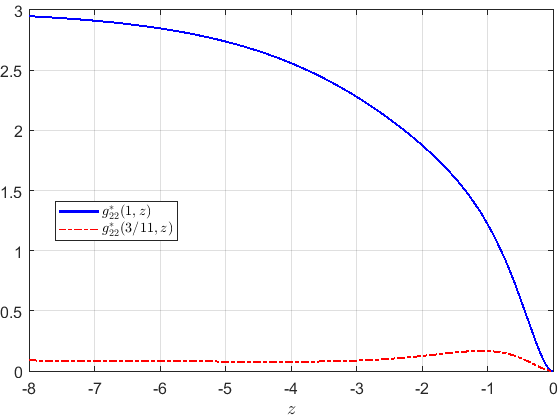}}
	\caption{Curves of comparison functions $g^*_{21}$ and $g^*_{22}$.}
	\label{fig: g21, g22 comparfuns}
\end{figure}

Proposition \ref{proposition: g21} shows that $\mathrm{Det}\kbrab{\mathcal{S}(D^{(2)};c_2,z)}>0$ for $c_2\in[\frac12,1]$ and $z<0$. Thus, the sequential principal minors of $\mathcal{S}(D^{(2)};c_2,z)$ are positive and then the differentiation matrix $D^{(2)}(c_2,z)$ is positive definite. 
Theorem \ref{thm: energy stability} gives the following result.
\begin{corollary}\label{corollary: ETDRK2} 
	The EERK2 method \eqref{scheme: EERK2 Butcher tableau} with $c_2\in[\frac12,1]$ preserves
	the energy dissipation law \eqref{problem: energy dissipation law} unconditionally
	at all stages in the sense that
	\begin{align*}		
		E[U^{n,j+1}]-E[U^{n,1}]\le&\,-\frac1{\tau}\sum_{k=1}^{j}\myinnerB{\delta_{\tau}U^{n,k+1},
			\sum_{\ell=1}^{k}d_{k\ell}^{(2)}(c_2,-{\tau}L_{\kappa})\delta_{\tau}U^{n,\ell+1}}\\
		&\,=-\frac1{\tau}\myinnerB{\delta_{\tau}\vec{U}_{n,j+1},
			D_{j}^{(2)}(c_2,-{\tau}L_{\kappa})\delta_{\tau}\vec{U}_{n,j+1}}
		\quad\text{for $1\le j\le 2$.}
	\end{align*}
\end{corollary}

\begin{figure}[htb!]
	\centering
	\subfigure[EERK2]{
		\includegraphics[width=2.15in]{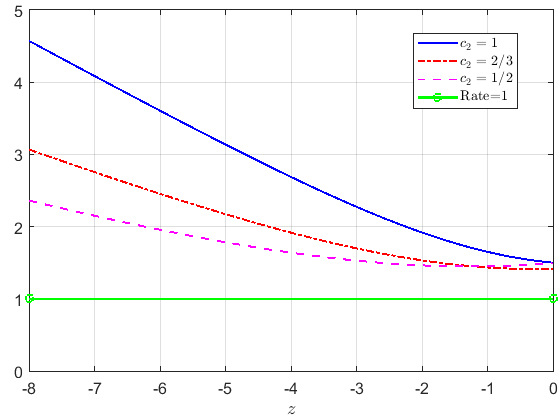}}\hspace{10mm}
	\subfigure[EERK2-w]{
		\includegraphics[width=2.15in]{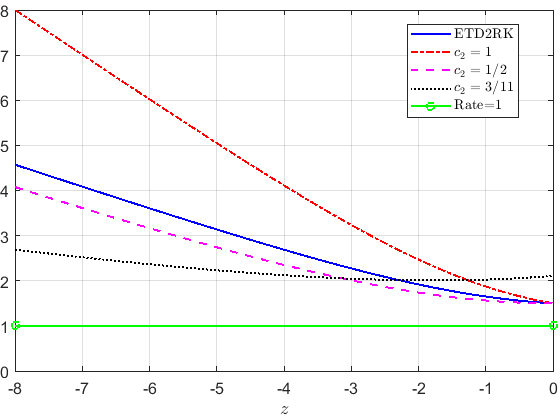}}	
	\caption{Averaged dissipation rates $\mathcal{R}^{(2)}(c_2,z)$ and $\mathcal{R}^{(2,w)}(c_2,z)$ for different abscissas $c_2$.}
	\label{fig: EERK2 average_rates}
\end{figure}

According to Lemma \ref{lemma: average dissipation rate}, 
the EERK2 method \eqref{scheme: EERK2 Butcher tableau} has the average dissipation rate
\begin{align}\label{def: EERK2 average rate}
	\mathcal{R}^{(2)}(c_2,z):=\frac{z}{2}+\frac1{2c_{2}\varphi_{1}(c_2z)}+\frac{c_{2}}{2\varphi_{2}(z)}
	\quad\text{for $c_2\in(0,1]$ and $z\le0$.}
\end{align}
It is different for different choices of $c_2$. One has
$$\lim\limits_{z\rightarrow0}\mathcal{R}^{(2)}(c_2,z)=\frac{1}{2c_2}+c_2
\quad\text{and}\quad
\lim\limits_{z\rightarrow-\infty}\mathcal{R}^{(2)}(c_2,z)=+\infty
\quad\text{for $c_2\in(0,1]$.}$$ 
For properly large time-step sizes,  
the ETD2RK method \eqref{scheme: ETDRK2} with the case $c_2=1$ has the largest dissipation rate, see Figure \ref{fig: EERK2 average_rates} (a), while the average dissipation rate of the case $c_2=\tfrac{1}{2}$ is much more closer to 1.
The former is preferred to preserve the unconditional contractivity, while the latter would be preferred to
preserve the energy dissipation law \eqref{problem: energy dissipation law} unconditionally. 
As seen, $\mathcal{R}^{(2)}(c_2,z)>1$ for $c_2\in[\tfrac{1}{2},1]$ and $z\le0$, so that 
the EERK2 method \eqref{scheme: EERK2 Butcher tableau} always generates a time ``ahead" in simulating the gradient system \eqref{problem: stabilized version}.

\subsection{Weak variants of EERK2 method}

For the EERK2-w methods with tableau \eqref{scheme: EERK2-weak Butcher tableau}, 
the associated Butcher-Diff tableau
\begin{align*}
	\text{EERK2-w Butcher-Diff:}\quad
	\begin{array}{c|cccccc}
		0 &  &   \\
		c_{2} & c_{2}\varphi_{1,2} &   \\
		\hline  & \brat{1-\frac{1}{2c_2}}\varphi_1-c_{2}\varphi_{1,2} & \frac{1}{2c_2}\varphi_1 
	\end{array}\quad.
\end{align*}
By the procedure \eqref{eq: orthogonal procedureII}, 
one can compute the associated DOC kernels
\begin{align*}
	\underline{\theta}_{11}^{(2,w)}=&\,\frac{1}{c_2\varphi_{1}(c_2z)},\quad 
	\underline{\theta}_{22}^{(2,w)}=\frac{2c_2}{\varphi_{1}(z)}\quad\text{and}\quad
	\underline{\theta}_{21}^{(2,w)}=\frac{2c_{2}}{\varphi_{1}(z)}
	+\frac{(1-2c_2)}{c_2\varphi_{1}(c_2z)}.	
\end{align*}
The definition \eqref{Def: Differential Matrix DII} gives 
the following one-parameter differentiation matrix
\begin{align*}
	D^{(2,w)}(c_2,z):
		=&\,\begin{pmatrix}
				\frac{1}{c_2\varphi_{1}(c_2z)}+\frac{z}{2} & 0  \\[4pt]	
				\frac{2c_{2}}{\varphi_{1}(z)}
				+\frac{(1-2c_2)}{c_2\varphi_{1}(c_2z)}
				+z & \frac{2c_2}{\varphi_{1}(z)}+\frac{z}{2}\\
			\end{pmatrix}.
\end{align*}	
It is not difficult to check that 
\begin{align*}
	\mathrm{Det}\kbrab{\mathcal{S}(D_1^{(2,w)};c_2,z)}=&\,\frac{z(e^{c_2z}+1)}{2(e^{c_2z}-1)}\ge \frac1{c_2}\quad\text{for $c_2\in(0,1]$ and $z\le0$.}
\end{align*}	
To handle the second leading principal minor, we need the following result.
\begin{proposition}\label{proposition: g22}
	For the abscissa $c_2\in[\tfrac3{11},1]$ and $z<0$, it holds that
	\begin{align}\label{def: g22}
		g_{22}(c_2,z):=&\,
		-4 c_2^2 (e^{c_2z}-e^{z})^2+4 c_2 (1-e^{z})(1-e^{c_2 z+z})-(1-e^{z})^2>0.
	\end{align}	
\end{proposition}
\begin{proof}
	We consider the following comparison function
	\begin{align*}
		g^*_{22}(c_2,z):=&-4 c_2^2 \brab{e^{\tfrac{3z}{11}}-e^{z}}^2+4 c_2 (1-e^{z})\brab{1-e^{\tfrac{24z}{11}}}-(1-e^{z})^2\quad\text{for $z<0$.}
	\end{align*}
	For  $c_2\in[\frac{3}{11},1]$ and $z<0$, it is obvious that 
	\begin{align*}
		g_{22}(c_2,z)-g^*_{22}(c_2,z)=&4 c_2\brab{e^{\tfrac{3z}{11}}-e^{c_2 z}}e^z\kbraB{c_2 \brab{e^{-\tfrac{8z}{11}}+e^{(c_2-1)z}-2}+(1-e^{z})}\ge0.		
		\end{align*}	
Note that, $g^*_{22}$ is a concave, quadratic polynomial with respect to $c_2$ due to $\partial_{c_2}^2g^*_{22}<0$.
	Reminding that $\lim\limits_{z\rightarrow-\infty}g_{22}^*(c_2,z)=4c_2-1>0$ for $c_2\in[\frac{3}{11},1]$, it is not difficult to check that, 
	cf. Figure \ref{fig: g21, g22 comparfuns}(b), 
	\begin{align*}
		g^*_{22}(1,z)>0\quad\text{and}\quad
		g^*_{22}(\tfrac{3}{11},z)>0\quad\text{for $z<0$.}
	\end{align*}
	They imply that $g_{22}^*(c_2,z)>0$ and then $g_{22}(c_2,z)\ge g_{22}^*(c_2,z)>0$ for $c_2\in[\tfrac3{11},1]$ and $z<0$.
\end{proof}
Reminding the auxiliary function \eqref{def: g22} and Proposition \ref{proposition: g22}, we know that the second leading principal minor of $\mathcal{S}(D^{(2,w)};c_2,z)$ is positive, that is,
\begin{align*}
	\mathrm{Det}\kbrab{\mathcal{S}\brab{D^{(2,w)};c_2,z}}
	=&\,\frac{z^2g_{22}(c_2,z)}{4 (e^z-1)^2(e^{c_2 z}-1)^2}>0\quad\text{for $c_2\in[\tfrac3{11},1]$ and $z<0$.}
\end{align*}
Thus the sequential principal minors of $\mathcal{S}(D^{(2,w)};c_2,z)$ are positive
and then the differentiation matrix $D^{(2,w)}(c_2,z)$ is positive definite. 
Theorem \ref{thm: energy stability} gives the following result.
\begin{corollary}\label{corollary: ETD2RK-weak} 
	The EERK2-w method \eqref{scheme: EERK2-weak Butcher tableau} with $c_2\in[\tfrac3{11},1]$ preserves
	the original energy dissipation law \eqref{problem: energy dissipation law} 
	at all stages in the sense that
	\begin{align*}		
		E[U^{n,j+1}]-E[U^{n,1}]\le-\frac1{\tau}\sum_{k=1}^{j}\myinnerB{\delta_{\tau}U^{n,k+1},
			\sum_{\ell=1}^{k}d_{k\ell}^{(2,w)}(c_2,-{\tau}L_{\kappa})\delta_{\tau}U^{n,\ell+1}}
		\quad\text{for $1\le j\le 2$.}
	\end{align*}
\end{corollary}

 It is worth mentioning that the choice $c_2\in[\frac12,1]$ 
 makes the EERK2-w method \eqref{scheme: EERK2-weak Butcher tableau} unconditionally contractive \cite{MasetZennaro:2009MCOM}.
Similar to the ETD2RK scheme \eqref{scheme: ETDRK2}, the EERK2-w method \eqref{scheme: EERK2-weak Butcher tableau}  arrives at the following weak variant of ETD2RK scheme for $c_2\in[\frac12,1]$,
\begin{subequations}\label{scheme: ETD2RK-weak}
	\begin{align}			
		U^{n,2}=&\,\varphi_{0}(-c_2{\tau}L_{\kappa})U^{n,1}
		+{\tau}c_2\varphi_{1}(-c_2{\tau}L_{\kappa})
		g_{\kappa}(U^{n,1}),\\
		U^{n,3}=&\,\varphi_{0}(-{\tau}L_{\kappa})U^{n,1}
		+(1-\frac{1}{2c_2}){\tau}\varphi_1(-{\tau}L_{\kappa})g_{\kappa}(U^{n,1})
		+\frac{{\tau}}{2c_2}\varphi_{1}(-{\tau}L_{\kappa})g_{\kappa}(U^{n,2}).
	\end{align}
\end{subequations}
As an advantage over the ETD2RK scheme \eqref{scheme: ETDRK2}, the weak variant \eqref{scheme: ETD2RK-weak} provides more choice of the abscissa $c_2$ to preserve both the contractivity and energy dissipation law unconditionally.

For the one-parameter EERK2-w method \eqref{scheme: EERK2-weak Butcher tableau}, the average dissipation rate
\begin{align}\label{def: EERK2-w average rate}
	\mathcal{R}^{(2,w)}(c_2,z):=\frac{z}{2}+\frac1{2c_{2}\varphi_{1}(c_2z)}+\frac{c_{2}}{\varphi_{1}(z)}
	\quad\text{for $c_2\in(0,1]$ and $z\le0$.}
\end{align}
It is easy to find that
$$\lim\limits_{z\rightarrow0}\mathcal{R}^{(2,w)}(c_2,z)=\frac{1}{2c_2}+c_2
\quad\text{and}\quad
\lim\limits_{z\rightarrow-\infty}\mathcal{R}^{(2,w)}(c_2,z)=+\infty
\quad\text{for $c_2\in(0,1]$.}$$ 
For properly large time-step sizes,  
the case $c_2=1$ has the largest dissipation rate, see Figure \ref{fig: EERK2 average_rates} (b), 
while the case $c_2=\tfrac{1}{2}$ has the smallest rate near $z=0$ and the case $c_2=\tfrac{3}{11}$ has the smallest rate for $z<-3$. More interestingly, 
the case $c_2=\tfrac{1}{2}$ seems superior to the ETD2RK method \eqref{scheme: ETDRK2} 
since the dissipation rate of the former is much closer to 1.
For all cases $c_2\in[\tfrac{3}{11},1]$, the average rate $\mathcal{R}^{(2,w)}(c_2,z)>1$ and 
the EERK2-w method \eqref{scheme: EERK2-weak Butcher tableau} always generates a time ``ahead" 
for the gradient flow system \eqref{problem: stabilized version}.

\subsection{Remarks for the three-stage EERK method}

We present some remarks for the one-parameter family 
of 3-stage EERK (called EERK2-S in short) method proposed by Strehmel and Weiner \cite{StrehmelWeiner:1992book} with second-order B-consistency,
\begin{align}\label{scheme: 3-stage EERK2}
	\begin{array}{c|cccccc}
		0 &  &   \\
		c_{2} & c_{2}\varphi_{1,2} &   \\
		1 & \varphi_{1,3}-\frac{1}{c_2}\varphi_{2,3} &  \frac{1}{c_2}\varphi_{2,3} \\
		\hline  & \varphi_1-\varphi_2 & 0 & \varphi_2
\end{array}\;.\end{align}
These methods satisfy all conditions up to stiff order two, see  (\ref{6.13a})-(\ref{6.13d}), described in Section \ref{section: third-order methods}.  By following the arguments in \cite{MasetZennaro:2009MCOM}, we know that the EERK2-S method \eqref{scheme: 3-stage EERK2} with the abscissa $c_2=1$ is also unconditionally contractive since $\varphi_1(z)\ge\varphi_2(z)$ for $z\le0$. 
 
 To establish the stage energy laws, we present the Butcher-Diff tableau
  \begin{align*}
 	\text{EERK2-S Butcher-Diff:}\quad
 	\begin{array}{c|cccccc}
 	0 &  &   \\
 	c_2 & c_2\varphi_{1,2} &   \\
 	1 & \varphi_{1}-\frac{1}{c_2}\varphi_{2}-c_{2}\varphi_{1,2}  &  \frac{1}{c_2}\varphi_{2} \\[5pt]
 	\hline \\[-10pt] &\frac{1-c_2}{c_2}\varphi_{2} & -\frac{1}{c_2}\varphi_{2} & \varphi_2
 \end{array}.
 \end{align*}
Note that the first three lines of this Butcher-Diff tableau are the same to the Butcher-Diff tableau \eqref{scheme: EERK2 Butcher-Diff} of the  EERK2 method \eqref{scheme: EERK2 Butcher tableau}. 
The stage energies $E[U^{n,j}]$ $(j=2,3)$ of the EERK2-S method \eqref{scheme: 3-stage EERK2} have the same dissipation rates to those of the two-stage EERK2 method. Thus, to preserve the energy dissipation law \eqref{problem: energy dissipation law} at all stages, the condition $c_2\in[\tfrac{1}{2},1]$ is also necessary for the EERK2-S method \eqref{scheme: 3-stage EERK2}.

\begin{figure}[htb!]
	\centering
	\subfigure[With contractivity]{
		\includegraphics[width=2.15in]{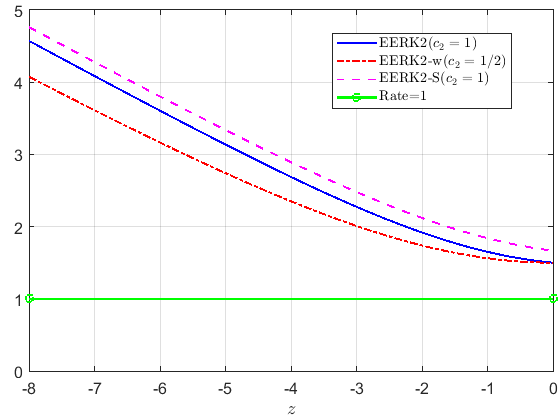}}\hspace{10mm}
	\subfigure[General case]{
		\includegraphics[width=2.15in]{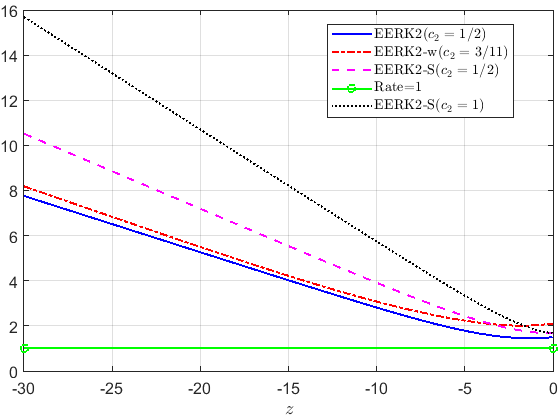}}
	\caption{Dissipation rate comparisons of EERK2, EERK2-w and EERK2-S methods.}
	\label{fig: comparison EERK2 average rates}
\end{figure}

However, the second-order EERK2-S method \eqref{scheme: 3-stage EERK2} would be not 
competitive for solving the gradient system \eqref{problem: stabilized version}.
Actually, one has the average dissipation rate
\begin{align*}
	\mathcal{R}^{(2,S)}(c_2,z)=\frac{z}2+\frac{1}{3c_{2}\varphi_{1}(c_2z)}
	+\frac{c_2}{3\varphi_{2}(z)}+\frac{1}{3\varphi_{2}(z)}\quad\text{for $z\le0$}.
\end{align*}
It is easy to obtain that 
$$\lim\limits_{z\rightarrow0}\mathcal{R}^{(2,S)}(c_2,z)=\frac{2}{3}+\frac2{3}\brab{c_2+\frac{1}{2c_2}}
\quad\text{and}\quad
\lim\limits_{z\rightarrow-\infty}\mathcal{R}^{(2,S)}(c_2,z)=+\infty
\quad\text{for $c_2\in(0,1]$.}$$ 
 In Figure \ref{fig: comparison EERK2 average rates}, we compare the 
average dissipation rates of the EERK2, EERK2-w and EERK2-S methods.  Taking into the contractivity account, 
we find that the EERK2-w method \eqref{scheme: EERK2-weak Butcher tableau} with $c_2=\tfrac{1}{2}$ generates the minimum time ``ahead" effect
since the average dissipation rate $\mathcal{R}^{(2,w)}(\frac12,z)$ has the smallest value for any time-step sizes, cf. Figure \ref{fig: comparison EERK2 average rates} (a),  while  the EERK2-S method with $c_2=1$ produces the maximum time ``ahead" effect. If the contractivity is not considered, Figure \ref{fig: comparison EERK2 average rates} (b) suggests that the EERK2 method \eqref{scheme: EERK2 Butcher tableau} with $c_2=\tfrac{1}{2}$ produces the minimum time ``ahead" effect among the three methods
since the average dissipation rate $\mathcal{R}^{(2)}(\frac12,z)$ has the smallest value.

Therefore, from the perspective of preserving the original energy dissipation rate,
the EERK2-S methods \eqref{scheme: 3-stage EERK2} with $c_2\in[\tfrac{1}{2},1]$ would be not 
competitive among second-order EERK methods for the gradient system \eqref{problem: stabilized version}.
As noted in \cite[subsection 5.2]{HochbruckOstermann:2005SINUM}, maybe more importantly, the three-stage the EERK2-S method \eqref{scheme: 3-stage EERK2} is not preferred for semilinear parabolic problems because they  are more computationally expensive than the two-stage methods in previous subsections.

\begin{table}[htb!]
	\begin{center}
		\caption{Choice of the abscissa $c_2$ in second-order EERK methods.
			\label{table: EERK2 Energy dissipation property}}
		\vspace*{0.3pt}
		\def\temptablewidth{1\textwidth}
		{\rule{\temptablewidth}{0.5pt}}
		\begin{tabular*}{\temptablewidth}{@{\extracolsep{\fill}}c|ccc}
			Method   & Contractivity & Energy law preserving & Best dissipation rate\\[4pt]	
			\hline					
			EERK2 \eqref{scheme: EERK2 Butcher tableau}   
			&$c_2=1$ 	  &$c_2\in[\tfrac{1}2,1]$&     $c_2=\tfrac{1}2$  \\[4pt]	
			\hline				
			EERK2-w \eqref{scheme: EERK2-weak Butcher tableau} 
			&$c_2\in[\tfrac{1}2,1]$ &$c_2\in[\tfrac{3}{11},1]$ & $c_2=\tfrac{3}{11}$ or $\tfrac{1}2$ 	  \\[4pt]			
			\hline		
			EERK2-S \eqref{scheme: 3-stage EERK2}  
			&$c_2=1$ & at least $c_2\in[\tfrac{1}2,1]$ &  $c_2=\tfrac{1}2$  	  \\	
		\end{tabular*}
		{\rule{\temptablewidth}{0.5pt}}
	\end{center}
\end{table}

As the end of this section, Table \ref{table: EERK2 Energy dissipation property} lists the abscissa choices  for the contractivity and the energy stability of three second-order EERK methods. Observations indicate that the abscissa condition maintaining the energy stability is often different from that preserving the contractivity, with the former being weaker than the latter. For the gradient flow system \eqref{problem: gradient flows}, certain time-stepping method without the contractivity does not necessarily violate the energy dissipation law \eqref{problem: energy dissipation law}. Actually, next subsection shows that
there are many third-order EERK methods preserving the energy dissipation law \eqref{problem: energy dissipation law}; however, as pointed out by \cite[Section 4]{MasetZennaro:2009MCOM}, scholars have not found that any EERK methods of stiff convergence order greater than two can maintain the contractivity.

\section{Discrete energy laws of third-order methods}
\label{section: third-order methods}
\setcounter{equation}{0}

Third-order methods require three internal stages, $s=3$, at least. The order conditions for three-stage methods \eqref{Scheme: general EERK3}	are given by \cite{HochbruckOstermann:2005SINUM}
\begin{subequations}
	\begin{align}
		a_{41}(-{\tau}L_{\kappa})+a_{42}(-{\tau}L_{\kappa})
		+a_{43}(-{\tau}L_{\kappa})&~=\varphi_1(-{\tau}L_{\kappa}),\label{6.13a}\\
		a_{42}(-{\tau}L_{\kappa})c_2+a_{43}(-{\tau}L_{\kappa})c_3
		&~=\varphi_2(-{\tau}L_{\kappa}),\label{6.13b}\\
		a_{21}(-{\tau}L_{\kappa})&~=c_2\varphi_1(-c_2{\tau}L_{\kappa}), \label{6.13c}\\
		a_{31}(-{\tau}L_{\kappa})+a_{32}(-{\tau}L_{\kappa})
		&~=c_3\varphi_1(-c_3{\tau}L_{\kappa}), \label{6.13d}\\
		a_{42}(-{\tau}L_{\kappa})c_2^2+a_{43}(-{\tau}L_{\kappa})c_3^2
		&~=2\varphi_3(-{\tau}L_{\kappa}), \label{6.13e}\\
		a_{42}(-{\tau}L_{\kappa})Jc_2^2\varphi_2(-c_2{\tau}L_{\kappa})
		+a_{43}(-{\tau}L_{\kappa})J\psi_{2,3}&~=0,\label{6.13f}
	\end{align}
\end{subequations}
where $J$ denotes arbitrary bounded operator and 
\begin{align*}
	\psi_{2,3}:=c_3^2\varphi_2(-c_3{\tau}L_{\kappa})-c_2a_{32}(-{\tau}L_{\kappa}).
\end{align*}
As pointed out in \cite[subsection 5.2]{HochbruckOstermann:2005SINUM},  
condition \eqref{6.13f} can be fulfilled by setting (I) $a_{42}=0$ and $\psi_{2,3}=0$; or (II)  $a_{42}=\gamma a_{43}$ and $c_2^2\varphi_2+\gamma \psi_{2,3}=0$.

The choice (I) leads to the following one-parameter family of method (called EERK3-1 in short)
\begin{align}\label{scheme: EERK3-HO1 Butcher}
	\begin{array}{c|cccccc}
			0 &  &   \\[2pt]
			c_{2} & c_{2}\varphi_{1,2} &   \\[3pt]
			\frac{2}{3} & \frac{2}{3}\varphi_{1,3}-\frac{4}{9c_2}\varphi_{2,3} &  \frac{4}{9c_2}\varphi_{2,3} \\[3pt]
			\hline  & \varphi_1-\frac{3}{2}\varphi_2 & 0 & \frac{3}{2}\varphi_2
		\end{array}\quad.
\end{align}
The other choice (II) leads to the two-parameter family of method (called EERK3-2) 
\begin{align}\label{scheme: EERK3-HO2 Butcher}
	\begin{array}{c|cccccc}
			0 &  &   \\[2pt]
			c_{2} & c_{2}\varphi_{1,2} &   \\[2pt]
			c_{3} & c_{3}\varphi_{1,3}-a_{32} &  \gamma c_2\varphi_{2,2}+\frac{c_3^2}{c_2}\varphi_{2,3} \\[3pt]
			\hline  & \varphi_1-a_{42}-a_{43} & \frac{\gamma}{\gamma c_2+c_3}\varphi_2 & \frac{1}{\gamma c_2+c_3}\varphi_2
		\end{array},
\end{align}
where the parameter $\gamma:=\frac{(3c_3-2)c_3}{(2-3c_2)c_2}$ for $c_2\neq\frac{2}{3}$ and $c_2\neq c_3$ (to ensure $a_{32}\neq0$). Also, it is to set $c_3\neq\tfrac{2}{3}$ since it degrades into the EERK3-1  method \eqref{scheme: EERK3-HO1 Butcher} if $c_3=\tfrac{2}{3}$.


In the literature, there are some related three-stage methods that involve the function $\varphi_3$. Cox and Matthews \cite{CoxMatthews:2002JCP} constructed the ETD3RK method with Butcher tableau
\begin{align}\label{scheme: ETD3RK Butcher}
	\begin{array}{c|cccccc}
		0 &  &   \\[2pt]
		\frac{1}{2} & \frac{1}{2}\varphi_{1,2} &   \\[2pt]
		1 & -\varphi_{1,3} &  2\varphi_{1,3} \\
		\hline  & 4\varphi_3-3\varphi_2+\varphi_1 & -8\varphi_3+4\varphi_2 & 4\varphi_3-\varphi_2
\end{array}\quad.\end{align}
This method satisfies the conditions \eqref{6.13a}-\eqref{6.13d}, while 
the conditions \eqref{6.13e}-\eqref{6.13f} are satisfied only in a very weak form (setting $L_{\kappa}=0$).
As a variant of the commutator-free Lie group CF3 method due to Celledoni, Marthinsen, and Owren \cite{CelledoniMarthinsenOwren:2003}, the so-called ETD2CF3 method  is given by 
\begin{align}\label{scheme: ETD2CF3 Butcher}
	\begin{array}{c|cccccc}
		0 &  &   \\[2pt]
		\frac{1}{3} & \frac{1}{3}\varphi_{1,2} &   \\[3pt]
		\frac{2}{3} & \frac{2}{3}\varphi_{1,3}-\frac{4}{3}\varphi_{2,3} &  \frac{4}{3}\varphi_{2,3} \\[3pt]
		\hline  & \varphi_1-\frac{9}{2}\varphi_2+9\varphi_3 & 6\varphi_2-18\varphi_3 & -\frac{3}{2}\varphi_2+9\varphi_3
\end{array}\quad.\end{align}
The  ETD2CF3 method satisfies the conditions \eqref{6.13a}-\eqref{6.13c}, while conditions \eqref{6.13d} and \eqref{6.13f} are satisfied in the weak form (setting $L_{\kappa}=0$).

\subsection{Simplified procedure and two simple cases}

As seen, the method coefficients (and the corresponding differentiation matrix as well) of high-order EERK methods are always more complex than those of lower order methods. Therefore certain symbolic computation system, such as the \emph{Wolfram Mathematica}, will be employed to assist our theoretical derivations for stage energy dissipation laws.  To do that, we summarize our theoretical framework in the subsections 2.1 and 2.2 as follows:
\begin{description}
\item[\textbf{Step1}.] Compute the differentiation  matrix $D(z)$ 
defined via \eqref{Def: Differential Matrix DII}, or
	$$D(z)=\bra{E_s^{-1}A(z)}^{-1}+zE_s-\frac{z}{2}I=A(z)^{-1}E_s+zE_s-\frac{z}{2}I,$$ where
	$A(z):=(\mathrm{a}_{ij})_{s\times s}$ is the coefficient matrix with $\mathrm{a}_{ij}:=a_{i+1,j}(z)$ for $1\le i,j\le s$ and $E_s:=(1_{i\ge j})_{s\times s}$ is the lower triangular matrix full of element 1. 
\item[\textbf{Step2}.] Compute the $j$-th leading principal minors $\mathrm{Det}\kbra{\mathcal{S}(D_j;z)}$ for $1\le j\le s$ and check the positive definiteness of the symmetric matrix $\mathcal{S}\brat{D;z}$ for $z\le0$. 
\item[\textbf{Step3}.] Establish the stage energy dissipation laws if $\mathcal{S}\brat{D;z}$ is positive definite and compute the average dissipation rate $\mathcal{R}(z)$ using the coefficients $a_{i+1,i}$ $(1\le i\le s)$ of the  EERK method. 
\end{description}

In the following, we will apply the procedure (\textbf{Step1})-(\textbf{Step3}) to examine the above third-order EERK methods and pick out those preserving the energy dissipation law \eqref{problem: energy dissipation law} unconditionally.

\begin{figure}[htb!]
	\centering
	\subfigure[ETD3RK \cite{CoxMatthews:2002JCP}]
	{\includegraphics[width=2.15in]{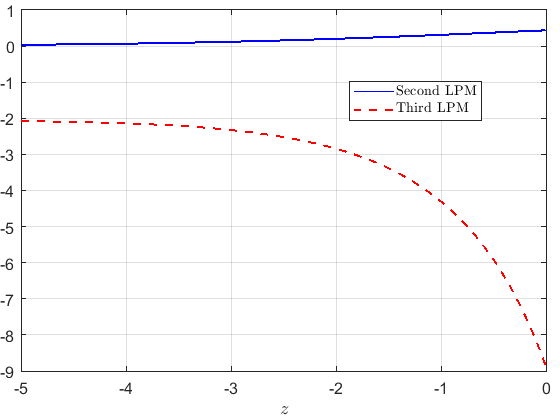}}\hspace{10mm}
	\subfigure[ETD2CF3 \cite{CelledoniMarthinsenOwren:2003}]
	{\includegraphics[width=2.15in]{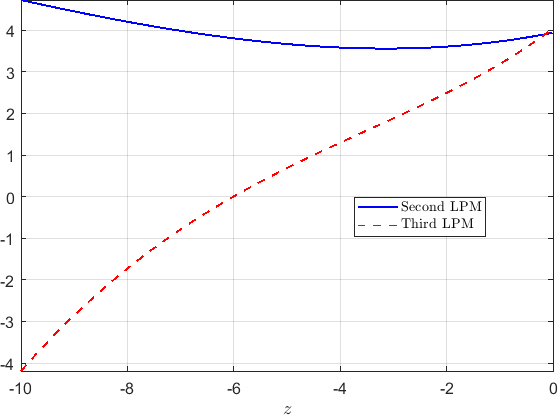}}
	\caption{Leading principal minors (LPM) of associated differential matrices.}
	\label{fig: principal minors ETD3RK and ETD2CF3}
\end{figure}

Before searching for some third-order EERK methods that preserve the energy dissipation law \eqref{problem: energy dissipation law} unconditionally, we first examine the well-known ETD3RK 
\eqref{scheme: ETD3RK Butcher} and ETD2CF3 \eqref{scheme: ETD2CF3 Butcher} 
methods by computing the associated differential matrices $D^{(3,e)}(z)$ and $D^{(3,f)}(z)$. Figure \ref{fig: principal minors ETD3RK and ETD2CF3} shows that the determinant of $\mathcal{S}(D^{(3,e)};z)$ is always negative for $z<0$ and the differential matrix $D^{(3,e)}(z)$ is not positive definite for any $z<0$. The determinant of $\mathcal{S}(D^{(3,f)};z)$ is always negative and the differential matrix $D^{(3,f)}(z)$ is not positive definite for $z<-6$. That is to say, when applied to the gradient system \eqref{problem: gradient flows}, 
both methods may destroy the energy dissipation law \eqref{problem: energy dissipation law} (especially for large time-step sizes) no matter how large the stabilization parameter $\kappa$ we set in \eqref{def: stabilized parameter}.

\subsection{One-parameter EERK3 methods}

For the $c_2$-parameterized EERK3-1 method \eqref{scheme: EERK3-HO1 Butcher},
one has the following differentiation matrix
\begin{align}\label{def: EERK3-1 differentiation matrix}
	D^{(3,1)}(c_2,z):=
	\begin{pmatrix}
		\frac{1}{c_2\varphi_1(c_2z)}+\frac{z}{2} & 0 & 0 \\
		\frac{9c_2}{4 \varphi_2(\frac{2 z}{3})}
		+\frac{1}{c_2\varphi_1(c_2z)}
		-\frac{3 \varphi_1(\frac{2 z}{3})}{2 \varphi_2(\frac{2 z}{3})\varphi_1(c_2z)}+z
		& \frac{9 c_2}{4 \varphi_2(\frac{2 z}{3})}+\frac{z}{2} & 0 \\
		\frac{2c_2\varphi_1(c_2z)-2 \varphi_1(z)+3 \varphi_2(z)}{3 c_2\varphi_1(c_2z) \varphi_2(z)}+z & \frac{2}{3 \varphi_2(z)}+z & \frac{2}{3 \varphi_2(z)}+\frac{z}{2} \\
	\end{pmatrix}.
\end{align}
It is not difficult to check that 
\begin{align*}
	\mathrm{Det}\kbrab{\mathcal{S}(D_1^{(3,1)};c_2,z)}=&\,\frac{z(e^{c_2z}+1)}{2(e^{c_2z}-1)}\ge \frac1{c_2}\quad\text{for $c_2\in(0,1]$ and $z\le0$.}
\end{align*}	
By using the auxiliary function $g_{31}$ in \eqref{def: g31} and Proposition \ref{proposition: g31}, 
one has
\begin{align*}
	\mathrm{Det}\kbrab{\mathcal{S}\brab{D_2^{(3,1)};c_2,z}}
	=&\,\frac{ (e^{c_2 z}-1)^{-2}z^2}{4(2 z-3 e^{\frac{2 z}{3}}+3)^2}g_{31}(c_2,c_2,z)>0\quad\text{ for  $c_2\in[\tfrac{4}{9},1]$ and $z<0$.}
\end{align*}
Also, the third leading principal minor of $\mathcal{S}(D^{(3,1)};c_2,z)$ is given by
\begin{align*}
	\mathrm{Det}\kbrab{\mathcal{S}\brab{D^{(3,1)};c_2,z}}
	=&\,\frac{z^4(e^{c_2 z}-1)^{-2}g_{32}(c_2,c_2,z)}{72(z-e^z+1)^2(2 z-3 e^{\frac{2 z}{3}}+3)^2}>0\quad\text{ for  $c_2\in[\tfrac{4}{9},1]$ and $z<0$,}
\end{align*}
where $g_{32}$ is defined by \eqref{def: g31} and Proposition \ref{proposition: g32} has been used.
Then the sequential principal minors of $\mathcal{S}(D^{(3,1)};c_2,z)$ are positive,
and then the differentiation matrix $D^{(3,1)}(c_2,z)$ is positive definite. 
Theorem \ref{thm: energy stability} gives the following result.
\begin{corollary}\label{corollary: EERK3-1} 
	The one-parameter EERK3-1 method \eqref{scheme: EERK3-HO1 Butcher} with $c_2\in[\tfrac{4}{9},1]$ preserves
	the energy dissipation law \eqref{problem: energy dissipation law} 
	at all stages in the sense that
	\begin{align*}		
		E[U^{n,j+1}]-E[U^{n,1}]\le-\frac1{\tau}\sum_{k=1}^{j}\myinnerB{\delta_{\tau}U^{n,k+1},
			\sum_{\ell=1}^{k}d_{k\ell}^{(3,1)}(c_2,-{\tau}L_{\kappa})\delta_{\tau}U^{n,\ell+1}}
		\quad\text{for $1\le j\le 3$.}
	\end{align*}
\end{corollary}

\begin{figure}[htb!]
	\centering	
	\subfigure[$\mathcal{R}^{(3,1)}(c_2,z)$]{
		\includegraphics[width=2.15in]{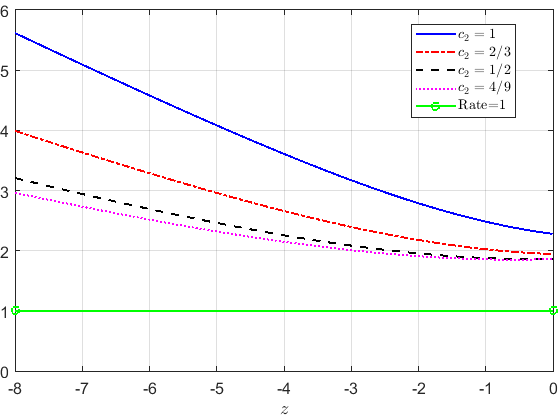}}\hspace{10mm}
	\subfigure[$\mathcal{R}^{(3,2)}(c_2,c_3,z)$]{
		\includegraphics[width=2.15in]{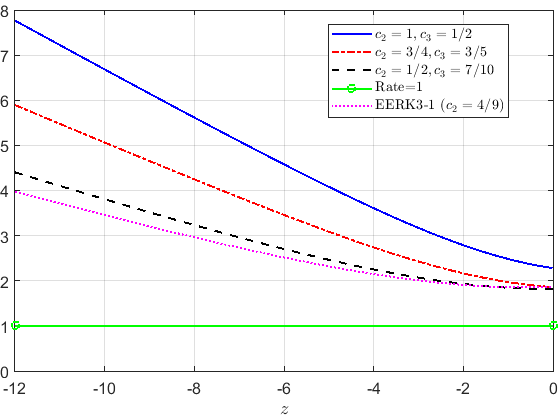}}
	\caption{Averaged dissipation rates of EERK3-1 and EERK3-2 methods.}
	\label{fig: dissipation rates of EERK3-HO}
\end{figure}

For the one-parameter EERK3-1 method \eqref{scheme: EERK3-HO1 Butcher}, the average dissipation rate
\begin{align*}
	\mathcal{R}^{(3,1)}(c_2,z)=\frac{z}{2}+\frac{1}{3c_{2}\varphi_{1}(c_2z)}
	+\frac{3c_2}{4\varphi_{2}(\tfrac{2z}{3})}+\frac{2}{9\varphi_{2}(z)}.
\end{align*}
It is easy to find that
$$\lim\limits_{z\rightarrow0}\mathcal{R}^{(3,1)}(c_2,z)=\frac{3c_2}{2}+\frac{1}{3c_2}+\frac{4}{9}
\quad\text{and}\quad
\lim\limits_{z\rightarrow-\infty}\mathcal{R}^{(3,1)}(c_2,z)=+\infty
\quad\text{for $c_2\in[\tfrac{4}{9},1]$.}$$ 
As seen in Figure \ref{fig: dissipation rates of EERK3-HO} (a), the case $c_2=1$ has the largest dissipation rate,  
while the case $c_2=\tfrac{4}{9}$ has the smallest rate. For all cases $c_2\in[\tfrac{4}{9},1]$, the average dissipation rate $\mathcal{R}^{(3,1)}(c_2,z)>1$ and 
the EERK3-1 method \eqref{scheme: EERK3-HO1 Butcher} always generates a time ``ahead" 
for the gradient system \eqref{problem: stabilized version}.

\subsection{Two-parameter EERK3 methods}
 
For the two-parameter EERK3-2 method \eqref{scheme: EERK3-HO2 Butcher}, it is reasonable to properly confine ourselves possible choices of the abscissas $c_2$ and $c_3$. According to Lemma \ref{lemma: average dissipation rate}, we consider the following average dissipation rate
\begin{align}\label{def: EERK3-HO2 dissipation rate}
	\mathcal{R}^{(3,2)}(c_2,c_3,z)=z+\frac{1}{3c_{2}\varphi_{1}(c_2z)}
	+\frac{1}{3\gamma c_2\varphi_{2}(c_2z)+\frac{3c_3^2}{c_2}\varphi_{2}(c_3z)}+\frac{\gamma c_2+c_3}{3\varphi_2(z)},
\end{align}
where $\gamma:=\frac{(3c_3-2)c_3}{(2-3c_2)c_2}$ for 
$c_2\neq\frac{2}{3}$, $c_2\neq c_3$ and $c_3\neq\frac{2}{3}$ (the choice $c_3:=\tfrac{2}{3}$ gives $\gamma=0$ and the method reduces to the EERK3-1  method \eqref{scheme: EERK3-HO1 Butcher}). It is not difficult to check that the following abscissa (necessary) condition 
\begin{align}\label{cond: EERK3-2 dissipation rate}
	\frac{6 c_3 (c_2-c_3)}{(3 c_2-2)}-1+\frac{2c_2(3 c_2-2)}{3 c_3 (c_2-c_3)}\ge0\quad\text{for 
		$c_2\neq\frac{2}{3}$ and $c_3\neq c_2$}
\end{align}
is sufficient to ensure that $\lim\limits_{z\rightarrow-\infty}\mathcal{R}^{(3,2)}(c_2,c_3,z)\ge0$.

In general, the quartic inequality \eqref{cond: EERK3-2 dissipation rate} gives rise to some theoretical trouble in proving the positive definiteness of the associated differentiation matrix $D^{(3,2)}(c_2,c_3,z)$. Actually, to prove the positivity of the second leading principal minor $\mathrm{Det}\kbrab{\mathcal{S}(D_2^{(3,2)};c_2,c_3,z)}$, one has to handle sixth degree polynomials with respect to $c_2$ or $c_3$ (while, in the previous subsection 3.1, only second degree polynomials, see Propositions \ref{proposition: g31} and \ref{proposition: g32}, should be handled to determine the positive definiteness of the differentiation matrix $D^{(3,1)}(c_2,z)$ of the EERK3-1 method); and the third leading principal minor $\mathrm{Det}\kbrab{\mathcal{S}(D_3^{(3,2)};c_2,c_3,z)}$  involves eighth degree polynomials with respect to $c_2$ or $c_3$. 

Since we are not able to present a complete discussion on the choices of $c_2$ and $c_3$ for the energy stability of EERK3-2 method \eqref{scheme: EERK3-HO2 Butcher}, this subsection examines some of concrete examples  falling into three cases listed as follows: (a) $c_2=1$ with $c_3\in(0,1]$, (b) $c_2\in(\tfrac{2}{3},1)$ with $c_3\in(0,c_2)$, and (c) $c_2\in(0,\tfrac{2}{3})$ with $c_3\in(c_2,1]$. The settings in all cases are necessary for the condition \eqref{cond: EERK3-2 dissipation rate}.

\begin{figure}[htb!]
	\centering
	\subfigure[$z_0=-1$]{
		\includegraphics[width=2.15in]{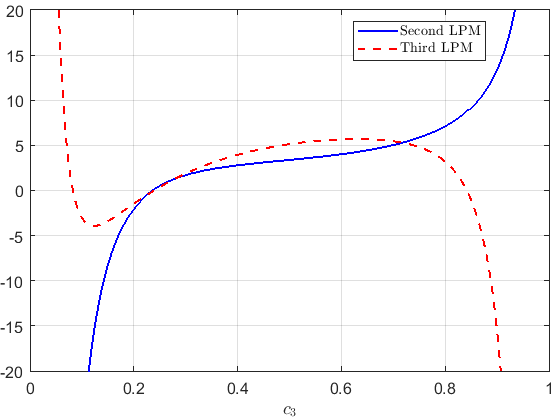}}\hspace{10mm}
	\subfigure[$z_0=-10$]{
		\includegraphics[width=2.15in]{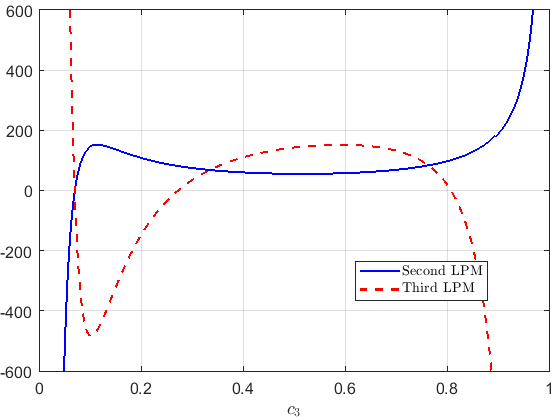}}
	\caption{Leading principal minors 
		$\mathrm{Det}\kbrat{\mathcal{S}(D_j^{(3,2)};1,c_3,z_0)}$ for $j=2,3$.}
	\label{fig: minors EERK3_HO-II_c2_1}
\end{figure}

\subsubsection{The case $c_2=1$ with $c_3\in(0,1]$} 
At first glance, the condition  \eqref{cond: EERK3-2 dissipation rate} always holds if $c_2=1$.
	We are to choose $c_3=\tfrac{1}{2}$ according to some numerical tests, cf. Figure \ref{fig: minors EERK3_HO-II_c2_1} (a)-(b), where the second and third leading principal minors of $\mathcal{S}(D_2^{(3,2)};1,c_3,z)$ are depicted for $z_0=-1$ and $z_0=-10$, respectively. 
	
	It is not difficult to check that 
	\begin{align*}
		\mathrm{Det}\kbrab{\mathcal{S}(D_1^{(3,2)};1,\tfrac{1}{2},z)}=&\,\frac{z(e^{z/2}+1)}{2(e^{z/2}-1)}>0 \quad\text{for $z\le0$.}
	\end{align*}	
	Reminding the auxiliary function  $g_{41}$ and $g_{42}$ defined by \eqref{def: g41}-\eqref{def: g42}, we know that the second and third leading principal minors of $\mathcal{S}(D^{(3,2)};1,\tfrac{1}{2},z)$ are positive, that is,
	\begin{align*}
		\mathrm{Det}\kbrab{\mathcal{S}\brab{D_2^{(3,2)};1,\tfrac{1}{2},z}}
		=&\,\frac{(e^{z/2}-1)^{-2}(e^{z/2}+1)^{-2}z^2}{4 (-3 z+4 e^{z/2}+e^z-5)^{2}}g_{41}(z)>0\quad\text{ for $z<0$,}\\
		\mathrm{Det}\kbrab{\mathcal{S}\brab{D^{(3,2)};1,\tfrac{1}{2},z}}
		=&\,\frac{ (e^z-1)^{-2}(z-e^z+1)^{-2}z^4}{128(-3 z+4 e^{z/2}+e^z-5)^{2}}g_{42}(z)>0\quad\text{ for $z<0$,}
	\end{align*}
	where Proposition \ref{proposition: g41 g42} has been used.
	Then the sequential principal minors of $\mathcal{S}(D^{(3,2)};1,\tfrac{1}{2},z)$ are positive
	and then the differentiation matrix $D^{(3,2)}(1,\tfrac{1}{2},z)$ is positive definite. 
	Theorem \ref{thm: energy stability} gives the following result.
	\begin{corollary}\label{corollary: EERK3-2a} 
		The two-parameter EERK3-2 method \eqref{scheme: EERK3-HO2 Butcher} with $c_2=1$ and $c_3=\tfrac{1}{2}$ preserves
		the energy dissipation law \eqref{problem: energy dissipation law} 
		at all stages in the sense that
		\begin{align*}		
			E[U^{n,j+1}]-E[U^{n,1}]\le-\frac1{\tau}\sum_{k=1}^{j}\myinnerB{\delta_{\tau}U^{n,k+1},
				\sum_{\ell=1}^{k}d_{k\ell}^{(3,2)}(1,\tfrac{1}{2},-{\tau}L_{\kappa})\delta_{\tau}U^{n,\ell+1}}
			\quad\text{for $1\le j\le 3$.}
		\end{align*}
	\end{corollary}

\subsubsection{The case $c_2\in(\tfrac{2}{3},1)$ with $c_3\in(0,c_2)$} 
\begin{figure}[htb!]
	\centering
	\subfigure[$z_0=-20$]{
		\includegraphics[width=2.15in]{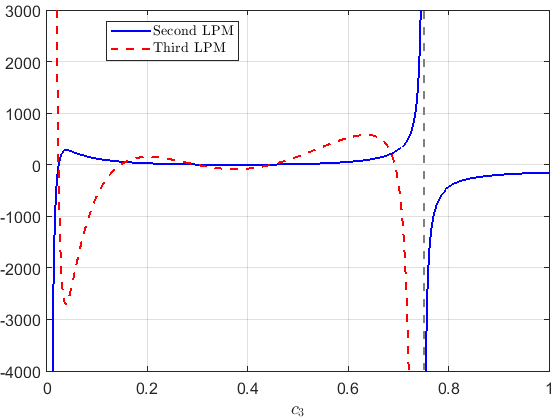}}\hspace{10mm}
	\subfigure[$z_0=-30$]{
		\includegraphics[width=2.15in]{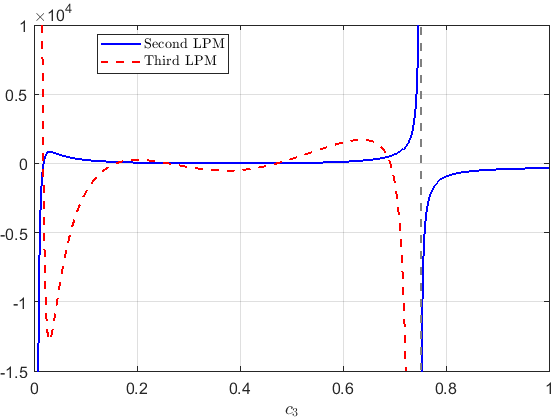}}
	\caption{Leading principal minors 
		$\mathrm{Det}\kbrat{\mathcal{S}(D_j^{(3,2)};\tfrac{3}4,c_3,z_0)}$ for $j=2,3$.}
	\label{fig: minors EERK3_HO-II_c2_075}
\end{figure}

We consider $c_2=\tfrac{3}4$ and choose $c_3=\tfrac{3}5$ according to some numerical tests, cf. Figure \ref{fig: minors EERK3_HO-II_c2_075} (a)-(b), where the second and third leading principal minors of $\mathcal{S}(D_2^{(3,2)};\tfrac{3}4,c_3,z)$ are depicted for $z_0=-20$ and $z_0=-30$, respectively. 

It is not difficult to check that 
\begin{align*}
	\mathrm{Det}\kbrab{\mathcal{S}(D_1^{(3,2)};\tfrac{3}4,\tfrac{3}5,z)}=&\,\frac{z(e^{3z/4}+1)}{2(e^{3z/4}-1)}>0 \quad\text{for $z\le0$.}
\end{align*}	
Reminding the auxiliary function  $g_{51}$ and $g_{52}$ defined by \eqref{def: g51}-\eqref{def: g52}, we know that the second and third leading principal minors of $\mathcal{S}(D^{(3,2)};\tfrac{3}4,\tfrac{3}5,z)$ are positive, that is,
\begin{align*}
	\mathrm{Det}\kbrab{\mathcal{S}\brab{D_2^{(3,2)};\tfrac{3}4,\tfrac{3}5,z}}
	=&\,\frac{100(e^{\frac{3 z}{4}}-1)^{-2}z^2}{(27 z-16 e^{\frac{3 z}{4}}-25 e^{\frac{3 z}{5}}+41)^2}g_{51}(z)>0\quad\text{ for $z<0$,}\\
	\mathrm{Det}\kbrab{\mathcal{S}\brab{D^{(3,2)};\tfrac{3}4,\tfrac{3}5,z}}
	=&\,\frac{300(e^{\frac{3 z}{4}}-1)^{-2}(z-e^z+1)^{-2}z^4}{(27 z-16 e^{\frac{3 z}{4}}-25 e^{\frac{3 z}{5}}+41)^2}g_{52}(z)>0\quad\text{ for $z<0$,}
\end{align*}
where Proposition \ref{proposition: g51 g52} has been used.
Then the sequential principal minors of $\mathcal{S}(D^{(3,2)};\tfrac{3}4,\tfrac{3}5,z)$ are positive
and then the differentiation matrix $D^{(3,2)}(\tfrac{3}4,\tfrac{3}5,z)$ is positive definite. 
Theorem \ref{thm: energy stability} gives the following result.
\begin{corollary}\label{corollary: EERK3-2b} 
	The two-parameter EERK3-2 method \eqref{scheme: EERK3-HO2 Butcher} 
	with $c_2=\tfrac{3}4$ and $c_3=\tfrac{3}5$ preserves
	the energy dissipation law \eqref{problem: energy dissipation law} 
	at all stages in the sense that
	\begin{align*}		
		E[U^{n,j+1}]-E[U^{n,1}]\le-\frac1{\tau}\sum_{k=1}^{j}\myinnerB{\delta_{\tau}U^{n,k+1},
			\sum_{\ell=1}^{k}d_{k\ell}^{(3,2)}(\tfrac{3}4,\tfrac{3}5,-{\tau}L_{\kappa})\delta_{\tau}U^{n,\ell+1}}
		\quad\text{for $1\le j\le 3$.}
	\end{align*}
\end{corollary}

\subsubsection{The case $c_2\in(0,\tfrac{2}{3})$ with $c_3\in(c_2,1]$} 
	
		\begin{figure}[htb!]
		\centering
		\subfigure[$z_0=-1$]{
			\includegraphics[width=2.15in]{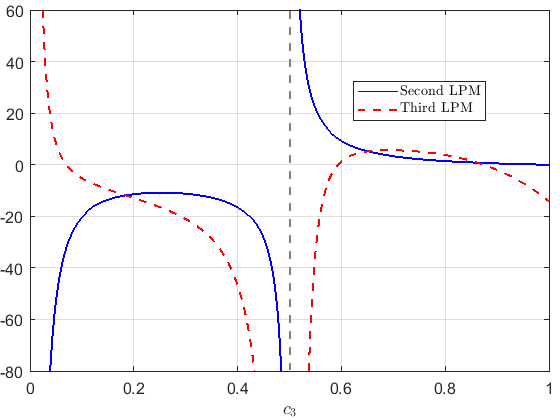}}\hspace{10mm}
		\subfigure[$z_0=-20$]{
			\includegraphics[width=2.15in]{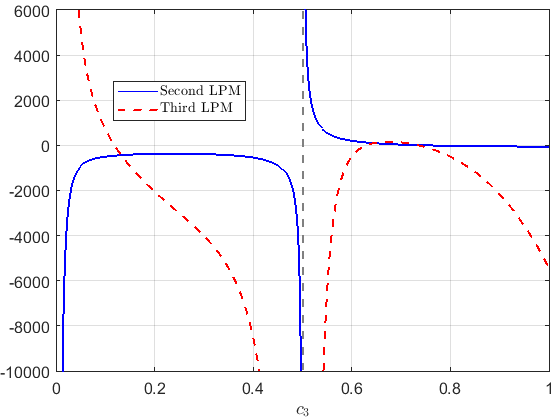}}
		\caption{Leading principal minors 
			$\mathrm{Det}\kbrat{\mathcal{S}(D_j^{(3,2)};\tfrac{1}2,c_3,z_0)}$ for $j=2,3$.}
		\label{fig: minors EERK3_HO-II_c2_05}
	\end{figure}
	
	We consider $c_2=\tfrac{1}2$ and choose the abscissa $c_3=\tfrac{7}{10}$ according to some numerical tests, cf. Figure \ref{fig: minors EERK3_HO-II_c2_05} (a)-(b), where the second and third leading principal minors of $\mathcal{S}(D_2^{(3,2)};\tfrac{1}2,c_3,z)$ are depicted for $z_0=-1$ and $z_0=-20$, respectively. 
	
	It is not difficult to check that 
	\begin{align*}
		\mathrm{Det}\kbrab{\mathcal{S}(D_1^{(3,2)};\tfrac{1}2,\tfrac{7}{10},z)}
		=&\,\frac{z(e^{z/2}+1)}{2(e^{z/2}-1)}>0 \quad\text{for $z\le0$.}
	\end{align*}	
	Reminding the auxiliary function  $g_{61}$ and $g_{62}$ defined by \eqref{def: g61}-\eqref{def: g62}, we know that the second and third leading principal minors of $\mathcal{S}(D^{(3,2)};\tfrac{1}2,\tfrac{7}{10},z)$ are positive, that is,
	\begin{align*}
		\mathrm{Det}\kbrab{\mathcal{S}\brat{D_2^{(3,2)};\tfrac{1}2,\tfrac{7}{10},z}}
		=&\,\frac{(e^{z/2}-1)^{-2}z^2}{(21 z-7 e^{z/2}-25 e^{\frac{7 z}{10}}+32)^2}g_{61}(z)>0\quad\text{ for $z<0$,}\\
		\mathrm{Det}\kbrab{\mathcal{S}\brat{D^{(3,2)};\tfrac{1}2,\tfrac{7}{10},z}}
		=&\,\frac{100(e^{z/2}-1)^{-2}(z-e^z+1)^{-2}z^4}{(21 z-7 e^{z/2}-25 e^{\frac{7 z}{10}}+32)^2}g_{62}(z)>0\quad\text{ for $z<0$,}
	\end{align*}
	where Proposition \ref{proposition: g61 g62} has been used.
	Then the sequential principal minors of $\mathcal{S}(D^{(3,2)};\tfrac{1}2,\tfrac{7}{10},z)$ are positive
	and then the differentiation matrix $D^{(3,2)}(\tfrac{1}2,\tfrac{7}{10},z)$ is positive definite. 
	Theorem \ref{thm: energy stability} gives the following result.
	\begin{corollary}\label{corollary: EERK3-2c} 
		The two-parameter EERK3-2 method \eqref{scheme: EERK3-HO2 Butcher} 
		with $c_2=\tfrac{1}2$ and $c_3=\tfrac{7}{10}$ preserves
		the energy dissipation law \eqref{problem: energy dissipation law} 
		at all stages in the sense that
		\begin{align*}		
			E[U^{n,j+1}]-E[U^{n,1}]\le-\frac1{\tau}\sum_{k=1}^{j}\myinnerB{\delta_{\tau}U^{n,k+1},
				\sum_{\ell=1}^{k}d_{k\ell}^{(3,2)}(\tfrac{1}2,\tfrac{7}{10},-{\tau}L_{\kappa})\delta_{\tau}U^{n,\ell+1}}
			\quad\text{for $1\le j\le 3$.}
		\end{align*}
	\end{corollary}

By the formula \eqref{def: EERK3-HO2 dissipation rate}, we compute the average dissipation rates $\mathcal{R}^{(3,2)}(c_2,c_3,z)$ for the above three examples in Corollaries \ref{corollary: EERK3-2a}, \ref{corollary: EERK3-2b} and \ref{corollary: EERK3-2c}.
For the above three cases of EERK3-2 method \eqref{scheme: EERK3-HO2 Butcher}, see Figure \ref{fig: dissipation rates of EERK3-HO} (b), the case $c_2=1$ 
with $c_3=\tfrac{1}{2}$ has the largest dissipation rate,  
while the case $c_2=\tfrac{1}2$ with $c_3=\tfrac{7}{10}$ has the smallest dissipation rate. Also, Figure \ref{fig: comparison EERK2 average rates} (b) suggests that the EERK3-1 method \eqref{scheme: EERK2 Butcher tableau} with $c_2=\tfrac{4}{9}$ produces the minimum time ``ahead" effect among the third-order EERK methods in this section
since the average dissipation rate $\mathcal{R}^{(3,1)}(\frac4{9},z)$ has the smallest value.
As the end of this section, Table \ref{table: EERK3 Energy dissipation property} 
summarizes the abscissa choices  for the energy stability of some third-order EERK methods. 

\begin{table}[htb!]	
	\centering
	\caption{The parameter choices in third-order EERK methods.
		\label{table: EERK3 Energy dissipation property}}			
	\vspace*{0.3pt}
	\def\temptablewidth{1\textwidth}
	{\rule{\temptablewidth}{0.5pt}}	
	\begin{threeparttable}
		\begin{tabular*}{\temptablewidth}{@{\extracolsep{\fill}}c|ccc}
			Method   & Unconditional energy law preserving & Best dissipation rate\\[6pt]
			\hline		
			EERK3-1 \eqref{scheme: EERK3-HO1 Butcher}  & $c_2\in[\tfrac{4}9,1]$ & $c_2=\tfrac{4}9$	\\[6pt]
			\hline
			\multirow{3}{*}{EERK3-2 \eqref{scheme: EERK3-HO2 Butcher}}    
			&  $c_2=1$, $c_3=\tfrac{1}{2}$  & 		
			\multirow{3}{*}{$c_2=\tfrac{1}2$, $c_3=\tfrac{7}{10}$} \\[4pt]		
			\cline{2-2}
						& $c_2=\tfrac{3}4$, $c_3=\tfrac{3}{5}$  & \\[6pt]	
			\cline{2-2}
			 			 & $c_2=\tfrac{1}2$, $c_3=\tfrac{7}{10}$  & 	\\[6pt]
			\hline	
		ETD3RK \eqref{scheme: ETD3RK Butcher}         &NPD*&    -- \\[6pt]
			\hline				
			ETD2CF3 \eqref{scheme: ETD2CF3 Butcher} 	 &NPD& --\\
		\end{tabular*}		
		{\rule{\temptablewidth}{0.5pt}}
		\footnotesize
		\tnote{* NPD means that there exists a $z_0<0$ such that the associated differential matrix $D(z_0)$ is not positive semi-definite.}	
	\end{threeparttable}		
\end{table}

{\color{blue}
\section{Numerical experiments}
For the sake of generality, we use the Cahn-Hilliard model,  $\partial_tu+\epsilon^2\Delta^2u=\Delta(u^3-u)$, to perform some numerical tests because this paper mainly focuses on the original energy dissipation properties of various EERK methods \eqref{Scheme: general EERK3}. Let $L_h$ be the discrete matrix of the Laplacian operator $-\Delta$. In such situation, the results of Theorem \ref{thm: energy stability} are valid by setting $L_\kappa:=\epsilon^2L_h^2+\kappa L_h$
and replacing the $L^2$ inner product $\myinnert{u,v}$ by the $H^{-1}$ inner product $\myinnert{u,v}_{-1}:=\myinnert{u,L_h^{-1}v}$ on the zero-mean function space $\{v|\myinnert{v,1}=0\}$. That is, we have the following discrete energy law
\begin{align*}			
	E[U^{n,j+1}]-E[U^{n,1}]\le&\,-\frac1{\tau}\sum_{k=1}^{j}\myinnerB{\delta_{\tau}U^{n,k+1},
		\sum_{\ell=1}^{k}d_{k\ell}\brab{-\epsilon^2\tau L_h^2-\kappa\tau L_h}\delta_{\tau}U^{n,\ell+1}}_{-1}
	\quad\text{for $1\le j\le s$.}
\end{align*}
In the first example, we examine the convergence by considering the EERK2-w methods \eqref{scheme: EERK2-weak Butcher tableau} and the EERK3-1 methods \eqref{scheme: EERK3-HO1 Butcher} for different abscissas $c_{2}$. In the second example, the energy dissipation rates are demonstrated for different choices of the method parameters including the  abscissa $c_{2}$, the stabilized parameter $\kappa$ and the time-step size $\tau$.


\subsection{Convergence tests}


\begin{example}\label{ex: FuShenYang} Consider the Cahn-Hilliard model on $\Omega=(0,2\pi)$ with $\epsilon=0.2$ subject to the initial data $u_0=0.5\sin(x)$ and Dirichlet boundary condition.
Always, we use the center difference approximation with the spacing $h=\pi/320$ for spatial discretization.
\end{example}


\begin{table}[htb!]	
	\centering
	\caption{EERK2-w errors with different abscissas $c_{2}$ for Example \ref{ex: FuShenYang}.}
	\label{tab: EERK2-w convergence rate}	
	\vspace*{0.3pt}
	\def\temptablewidth{1\textwidth}
	{\rule{\temptablewidth}{0.5pt}}	
	\begin{threeparttable}
		\begin{tabular*}{\temptablewidth}{@{\extracolsep{\fill}}c|cc|cc|cc|cc}
			\multirow{2}{*}{$\tau=0.01$} & \multicolumn{2}{c|}{$c_{2}=1$} & \multicolumn{2}{c|}{$c_{2}=\frac{3}{4}$} & \multicolumn{2}{c|}{$c_{2}=\frac{1}{2}$} & \multicolumn{2}{c}{$c_{2}=\frac{3}{11}$} \\[6pt]
			\cline{2-9}
			& $e(\tau)$ & Order & $e(\tau)$ & Order & $e(\tau)$ & Order & $e(\tau)$ & Order \\[6pt]
			\hline		
			$\tau$ & 6.106e-03 & - & 5.149e-03 & - & 4.122e-03 & - & 3.119e-03 & - \\[6pt]
			$\tau/2$ & 1.750e-03 & 1.80 & 1.462e-03 & 1.82 & 1.161e-03 & 1.83 & 8.756e-04 & 1.83 \\[6pt]
			$\tau/4$ & 4.744e-04 & 1.88 & 3.932e-04 & 1.89 & 3.098e-04 & 1.91 & 2.323e-04 & 1.91 \\[6pt]
			$\tau/8$ & 1.220e-04 & 1.96 & 1.001e-04 & 1.97 & 7.798e-05 & 1.99 & 5.756e-05 & 2.01 \\[6pt]
		\end{tabular*}		
		{\rule{\temptablewidth}{0.5pt}}
	\end{threeparttable}		
\end{table}

We always choose the final time $T=8$ and the stabilized parameter $\kappa=2$ in our tests. We run the EERK2-w methods \eqref{scheme: EERK2-weak Butcher tableau} with four different abscissas $c_2=1,\, \frac{3}{4}, \, \frac{1}{2}$ and $\frac{3}{11}$ for a small time step $\tau=0.001$. The four schemes work well and the corresponding solution and energy curves (omitted here) are hard to distinguish from each other. The numerical solution of the EERK2-w method with $c_2=\frac{3}{11}$ and $\tau=0.01/32$ is taken as the reference solution $u_h^{\star}$ in the convergence tests. The solution errors recorded in Table \ref{tab: EERK2-w convergence rate} are obtained on halving time steps $\tau=0.01/2^k$ for $k=0,1,\cdots,3$ and  the convergence order is computed by $\text{Order}\approx\log_2\bra{e(\tau)/e(\tau/2)}$
where $e(\tau)$ is the $L^\infty$ norm error defined by $e(\tau):=\max_{1\le n\le N}\mynorm{u_h^{n} - u_h^{\star}}_\infty$. The numerical results in Table \ref{tab: EERK2-w convergence rate} confirm the second-order time accuracy of the EERK2-w methods \eqref{scheme: EERK2-weak Butcher tableau}.

\begin{table}[htb!]	
	\centering
	\caption{EERK3-1 errors with different abscissas $c_{2}$ for Example \ref{ex: FuShenYang}.}
	\label{tab: EERK31 convergence rate}	
	\vspace*{0.3pt}
	\def\temptablewidth{1\textwidth}
	{\rule{\temptablewidth}{0.5pt}}	
	\begin{threeparttable}
		\begin{tabular*}{\temptablewidth}{@{\extracolsep{\fill}}c|cc|cc|cc|cc}
			\multirow{2}{*}{$\tau=0.01$} & \multicolumn{2}{c|}{$c_{2}=1$} & \multicolumn{2}{c|}{$c_{2}=\frac{2}{3}$} & \multicolumn{2}{c|}{$c_{2}=\frac{1}{2}$} & \multicolumn{2}{c}{$c_{2}=\frac{4}{9}$} \\[6pt]
			\cline{2-9}
			&$e(\tau)$ & Order & $e(\tau)$ & Order & $e(\tau)$ & Order & $e(\tau)$ & Order \\[6pt]
			\hline		
			$\tau$ & 6.369e-4 & - & 4.670e-4 & - & 3.708e-4 & - & 3.368e-4 & - \\[6pt]
			$\tau/2$ & 1.107e-4 & 2.52 & 7.922e-5 & 2.56 & 6.202e-5 & 2.58 & 5.604e-5 & 2.59 \\[6pt]
			$\tau/4$ & 1.737e-5 & 2.67 & 1.218e-5 & 2.70 & 9.425e-6 & 2.72 & 8.479e-6 & 2.72 \\[6pt]
			$\tau/8$ & 2.511e-6 & 2.79 & 1.729e-6 & 2.82 & 1.321e-6 & 2.83 & 1.183e-6 & 2.84 \\[6pt]
		\end{tabular*}
		{\rule{\temptablewidth}{0.5pt}}
	\end{threeparttable}		
\end{table}

It is interesting to note that, for each time step, the EERK2-w solutions with different $c_2$ have slight differences in precision, and the method with $c_2=\frac{3}{11}$ generates a bit more accurate solution than other cases. This interesting phenomenon is also observed from the solutions (omitted here) of EERK2 methods \eqref{scheme: EERK2 Butcher tableau}: the method with $c_2=\frac{1}{2}$ generates a bit more accurate solution than other choices $c_2>\frac{1}{2}$ including the widespread ETDRK method \cite{CoxMatthews:2002JCP} with $c_2=1$, at least for Example \ref{ex: FuShenYang}.  Coincidentally, the minimum abscissa choice preserving the energy dissipation law \eqref{problem: energy dissipation law} are $c_2=\frac{1}{2}$ for the EERK2 methods \eqref{scheme: EERK2 Butcher tableau} and $c_2=\frac{3}{11}$ for the EERK2-w methods \eqref{scheme: EERK2-weak Butcher tableau}, respectively, see Corollaries \ref{corollary: ETDRK2} and \ref{corollary: ETD2RK-weak}.

The numerical results of EERK3-1 methods \eqref{scheme: EERK3-HO1 Butcher}, listed in Table \ref{tab: EERK31 convergence rate}, are obtained in similar to those in Table \ref{tab: EERK2-w convergence rate}. As seen, the EERK3-1 methods always generate third-order solutions, at least for the smooth initial data. Also, we observe that the methods with $c_2>\frac{4}{9}$ generate a bit less accurate solution than the choice $c_2=\frac{4}{9}$,  which is the minimum abscissa preserving the energy dissipation law \eqref{problem: energy dissipation law} for the EERK3-1 methods, see Corollary \ref{corollary: EERK3-1}.

\subsection{Energy dissipation property}

\begin{example}\label{ex: Trefethen} Consider the Cahn-Hilliard model on $\Omega=(0,2\pi)$ with $\epsilon=0.2$ and zero-valued Dirichlet boundary condition subject to the following initial data,
	 see \cite[guide19]{DriscollHaleTrefethen:2014Chebfun}, 
\begin{align*}
	u_0=\frac13\tanh(2\sin x) - e^{-23.5(x-\frac{\pi}{2})^2}+e^{-27(x-4.2)^2} + e^{-38(x-5.4)^2}.
\end{align*}
We use the center difference approximation with the spacing $h=\pi/320$ for spatial discretization.
\end{example}

Taking the parameter $\kappa=2$ and the time-step $\tau=0.1$, we run the EERK2-w methods \eqref{scheme: EERK2-weak Butcher tableau} with four different abscissas $c_2=1,\frac{3}{4},\frac{1}{2},\frac{3}{11}$, and the EERK3-1 methods \eqref{scheme: EERK3-HO1 Butcher} with four different abscissas $c_2=1,\frac{2}{3},\frac{1}{2},\frac{4}{9}$ to the final time $T=160$. The corresponding numerical solution and discrete energy are depicted in Figures \ref{fig: EERK2w decay c2}-\ref{fig: EERK31 decay c2}, respectively. Taking a small step size $\tau=0.001$, we compute the reference solutions and energies (marked by ``Ref'' here and hereafter) by using the EERK2-w method with $c_2=\frac{3}{11}$ and the EERK3-1 method with $c_2=\frac{4}{9}$, respectively. As predicted by Corollaries \ref{corollary: ETD2RK-weak} and \ref{corollary: EERK3-1}, the original energies $E[u_h^n]$ generated by the two methods always decay over the time. 

As seen in Figure \ref{fig: EERK2w decay c2}(b), there are some obvious differences in energy dissipation rates for different abscissas $c_2$. It is not mysterious, according to the discrete energy law in Corollary \ref{corollary: ETD2RK-weak}, because the EERK2-w methods with different abscissas $c_2$ have different differentiation matrices $D^{(2,w)}(c_2,z)$. For this example, the discrete energy produced by the case $c_2=1$ decays fastest, while that generated by the case $c_2=\frac{3}{11}$ decays slowest. Qualitatively, they may be explained by the average dissipation rate $\mathcal{R}^{(2,w)}(c_2,z)$ in Figure \ref{fig: EERK2 average_rates}(b), in which we see that $\mathcal{R}^{(2,w)}(1,z)$ has the largest value and $\mathcal{R}^{(2,w)}(\tfrac{3}{11},z)$ has the smallest one for properly large $\abs{z}\ge4$. Similarly, the differences in energy dissipation rates for the EERK3-1 methods with different abscissas $c_2$, see Figure \ref{fig: EERK31 decay c2}(b), can be attributed to the differences of differentiation matrices $D^{(3,1)}(c_2,z)$ defined in \eqref{def: EERK3-1 differentiation matrix}. Also, they may be qualitatively explained by the average dissipation rate $\mathcal{R}^{(3,1)}(c_2,z)$ in Figure \ref{fig: dissipation rates of EERK3-HO}(a), in which we see that $\mathcal{R}^{(3,1)}(1,z)$ has the largest value and $\mathcal{R}^{(3,1)}(\tfrac{4}{9},z)$ has the smallest one for any $z<0$.

\begin{figure}[htb!]
  \centering
  \subfigure[final solution $u_h^N$]{
  \includegraphics[width=2.7in]{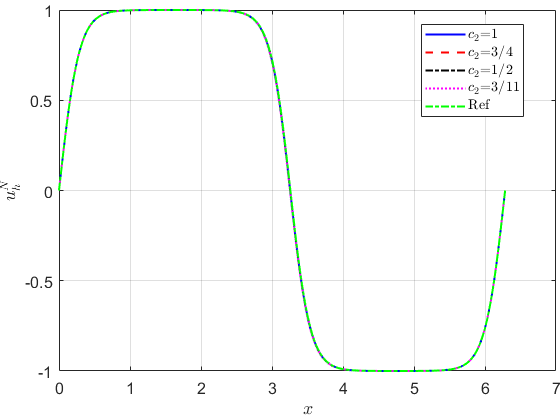}}\hspace{10mm}
  \subfigure[energy $E\kbrat{u_h^n}$]{
  \includegraphics[width=2.7in]{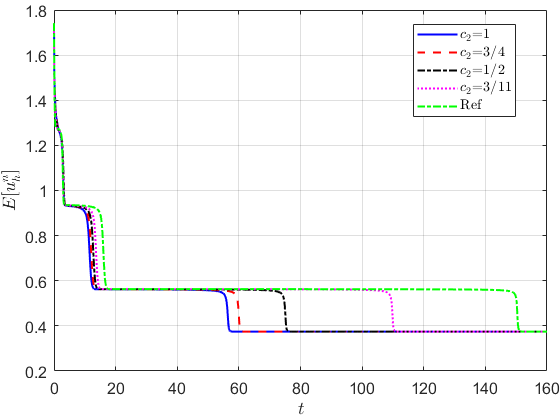}}
  \caption{Energy dissipation of the EERK2-w methods \eqref{scheme: EERK2-weak Butcher tableau} with $\kappa=2$ and $\tau=0.1$.}
  \label{fig: EERK2w decay c2}
\end{figure}
\begin{figure}[htb!]
  \centering
  \subfigure[final solution $u_h^N$]{
  \includegraphics[width=2.7in]{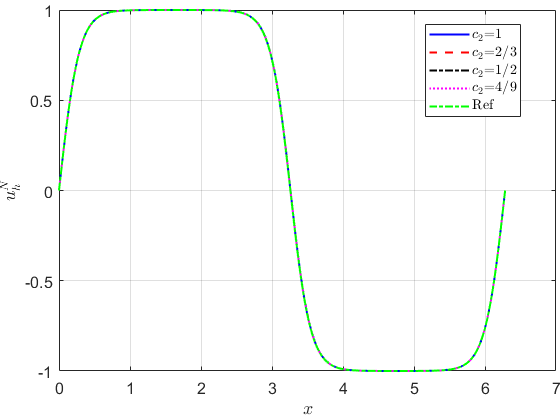}}\hspace{10mm}
  \subfigure[energy $E\kbrat{u_h^n}$]{
  \includegraphics[width=2.7in]{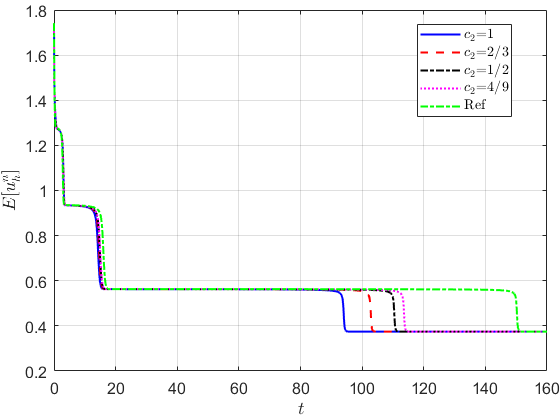}}
  \caption{Energy dissipation of the EERK3-1 methods \eqref{scheme: EERK3-HO1 Butcher} with $\kappa=2$ and $\tau=0.1$.}
  \label{fig: EERK31 decay c2}
\end{figure}


\begin{figure}[htb!]
	\centering
	\subfigure[final solution $u_h^N$]{
		\includegraphics[width=2.7in]{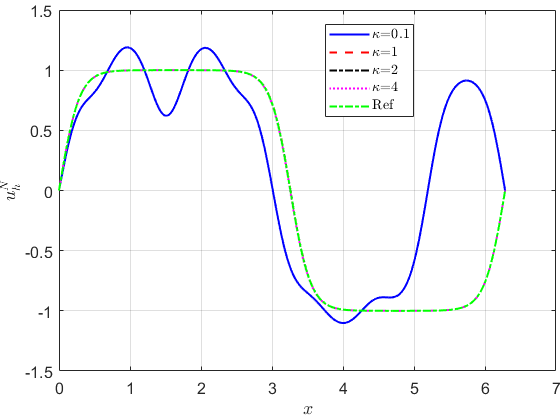}}\hspace{10mm}
	\subfigure[energy $E\kbrat{u_h^n}$]{
		\includegraphics[width=2.7in]{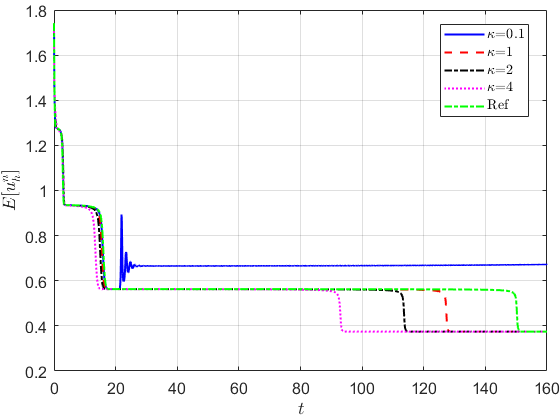}}
	\caption{Energy dissipation of the EERK3-1 method \eqref{scheme: EERK3-HO1 Butcher} with $c_2=\frac{4}{9}$ and $\tau=0.1$.}
	\label{fig: EERK31 decay kappa}
\end{figure}
\begin{figure}[htb!]
	\centering
	\subfigure[final solution $u_h^N$]{
		\includegraphics[width=2.7in]{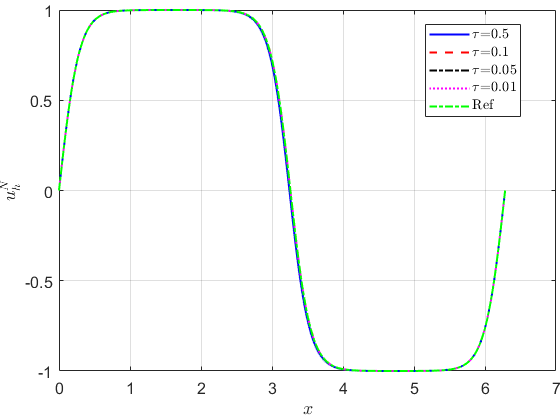}}\hspace{10mm}
	\subfigure[energy $E\kbrat{u_h^n}$]{
		\includegraphics[width=2.7in]{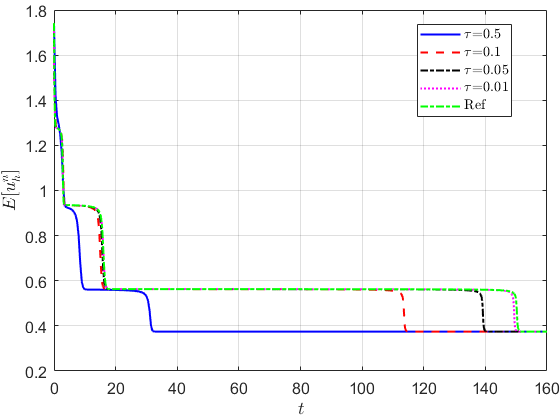}}
	\caption{Energy dissipation of the EERK3-1 method \eqref{scheme: EERK3-HO1 Butcher} with $c_2=\frac{4}{9}$ and $\kappa=2$.}
	\label{fig: EERK31 decay tau}
\end{figure}

Obviously, in addition to the different dissipation rates brought by the different choices of $c_2$, the time step size $\tau$ and stabilized parameter $\kappa$ also have some significant impacts on the discrete energy dissipation property. 
To explore the influence of stabilized parameter $\kappa$, we take a fixed time-step $\tau=0.1$ and run the EERK3-1 method with $c_2=\frac{4}{9}$ to the final time $T=160$ for four different parameters $\kappa=0.1,1,2$ and 4, cf. Figure \ref{fig: EERK31 decay kappa}, where the reference solution is computed with $\kappa=4$ and $\tau=0.001$. The discrete energy for $\kappa=0.1$ appears non-physical oscillations since the nonlinear stability could not be controlled by the small stabilized parameter. In practice, a properly large $\kappa$ is always necessary to maintain the stability especially when some large time step $\tau$ is employed.  With the increase of  $\kappa$, the energy curve appears some ``ahead" effect, that is, the discrete energy dissipates faster as the stabilization parameter $\kappa$ becomes larger.

Now we fix the stabilized parameter $\kappa=2$ and run the EERK3-1 method with $c_2=\frac{4}{9}$ for four different time steps $\tau=0.5, 0.1, 0.05$ and $0.01$, cf. Figure \ref{fig: EERK31 decay tau}, in which the reference solution is obtained with $\tau=0.001$. We see that, with the increase of time-step size, the energy curve shows some ``ahead" effect, that is, the discrete energy dissipate faster as the step size $\tau$ becomes larger. Note that, the numerical behaviors in Figure \ref{fig: EERK31 decay tau}(b) and \ref{fig: EERK31 decay kappa}(b) would be predictable by the average dissipation rate $\mathcal{R}^{(3,1)}(\tfrac{4}{9},z)$, see Figure \ref{fig: dissipation rates of EERK3-HO}(a), since it is increasing with respect to $\abs{z}$. Actually, we also run the EERK2-w method with $c_2=\frac12$ for Example \ref{ex: Trefethen} and find similar behaviors (omitted here) of the discrete energy curves for different time steps $\tau$ and different stabilized parameters $\kappa$. 


}

\section{Fourth-order EERK methods and concluding remarks}
\setcounter{equation}{0}

\begin{figure}[htb!]
	\centering
	\subfigure[Cox and Matthews (2002)]
	{\includegraphics[width=2.15in]{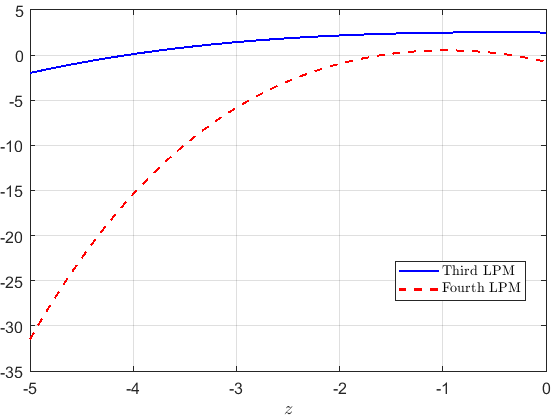}}\hspace{10mm}
	\subfigure[Krogstad (2005)]
	{\includegraphics[width=2.15in]{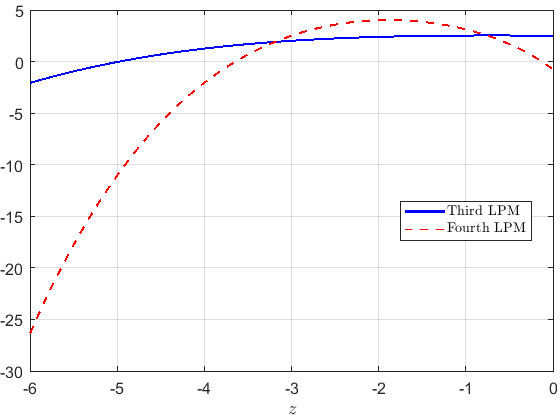}}\\
	\subfigure[Strehmel and Weine (1992)]
	{\includegraphics[width=2.15in]{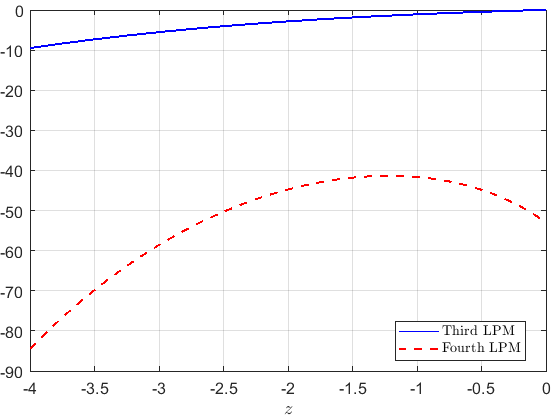}}\hspace{10mm}
	\subfigure[Hochbruck and Ostermann (2005)]
	{\includegraphics[width=2.15in]{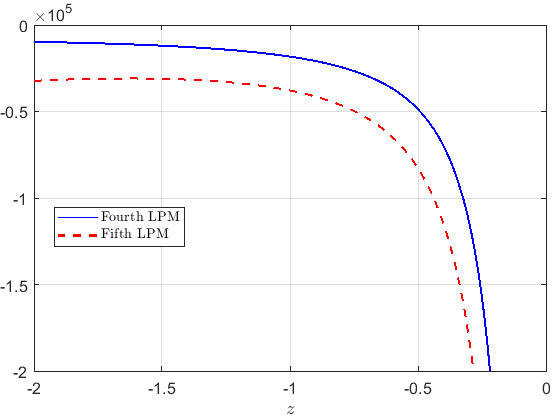}}
	\caption{Some leading principal minors (LPM) of associated differential matrices
			generated by existing fourth-order EERK methods in \cite{CoxMatthews:2002JCP,HochbruckOstermann:2005SINUM,StrehmelWeiner:1992book,Krogstad:2005JCP}.}
	\label{fig: principal minors of EERK4 methods}
\end{figure}

We consider firstly three four-stage fourth-order EERK methods from \cite{CoxMatthews:2002JCP,Krogstad:2005JCP,StrehmelWeiner:1992book}. As noted in \cite{HochbruckOstermann:2005SINUM}, these methods do not have the stiff order four although they show a
higher order of convergence (generically up to order four) under favorable circumstances.
The first one is the following exponential variant of the classical Runge-Kutta method developed by 
Cox and Matthews \cite{CoxMatthews:2002JCP}
\begin{align}\label{scheme: EERK4-Cox Butcher}
	\begin{array}{c|cccccc}
		0 &  &   \\[3pt]
		\frac{1}{2} & \frac{1}{2}\varphi_{1,2}     \\[3pt]
		\frac{1}{2} & 0 & \frac{1}{2}\varphi_{1,3} \\[3pt]
		1 & \frac{1}{2}\varphi_{1,3}(\varphi_{0,3}-1) & 0 & \varphi_{1,3} \\[3pt]
		\hline  & \varphi_1-3\varphi_2+4\varphi_3 & 2\varphi_2-4\varphi_3 
		& 2\varphi_2-4\varphi_3 & 4\varphi_3-\varphi_2
\end{array}\quad.\end{align}
The second one is the Krogstad's method \cite{Krogstad:2005JCP} given by
\begin{align}\label{scheme: EERK4-Krogstad Butcher}
	\begin{array}{c|cccccc}
		0 &  &   \\[3pt]
		\frac{1}{2} & \frac{1}{2}\varphi_{1,2}     \\[3pt]
		\frac{1}{2} & \frac{1}{2}\varphi_{1,3}-\varphi_{2,3} & \varphi_{2,3} \\[3pt]
		1 & \varphi_{1,4}-2\varphi_{2,4} & 0 & 2\varphi_{2,4} \\[3pt]
		\hline  & \varphi_1-3\varphi_2+4\varphi_3 & 2\varphi_2-4\varphi_3 & 2\varphi_2-4\varphi_3 & -\varphi_2+4\varphi_3
\end{array}\quad.\end{align}
The last is the following method from Strehmel and Weiner \cite[Example 4.5.5]{StrehmelWeiner:1992book}, \begin{align}\label{scheme: EERK4-Strehmel Butcher}
	\begin{array}{c|cccccc}
		0 &  &   \\[3pt]
		\frac{1}{2} & \frac{1}{2}\varphi_{1,2}     \\[3pt]
		\frac{1}{2} & \frac{1}{2}\varphi_{1,3}-\frac{1}{2}\varphi_{2,3} & \frac{1}{2}\varphi_{2,3} \\[3pt]
		1 & \varphi_{1,4}-2\varphi_{2,4} & -2\varphi_{2,4} & 4\varphi_{2,4} \\[3pt]
		\hline  & \varphi_1-3\varphi_2+4\varphi_3 & 0 & 4\varphi_2-8\varphi_3 & -\varphi_2+4\varphi_3
\end{array}\quad.
\end{align}
We compute the associated differential matrices $D^{(4,C)}(z)$, $D^{(4,K)}(z)$ and $D^{(4,S)}(z)$ of the above three methods \eqref{scheme: EERK4-Cox Butcher}-\eqref{scheme: EERK4-Strehmel Butcher}. Numerical results in Figure \ref{fig: principal minors of EERK4 methods} (a)-(c) show that the third and fourth leading principal minors are \lan{not always positive} for $z<0$. That is, the differential matrices $D^{(4,C)}(z)$, $D^{(4,K)}(z)$ and $D^{(4,S)}(z)$ are not positive (semi-)definite. It seems that these EERK methods would not be stabilized to preserve the energy dissipation law \eqref{problem: energy dissipation law} no matter how large the stabilization parameter $\kappa$ we set in \eqref{def: stabilized parameter}.

Hochbruck and Ostermann \cite{HochbruckOstermann:2005SINUM} constructed the following five-stage  EERK method which has been proved to have the stiff order four,
\begin{align}\label{eq: Butcher12}
	\begin{array}{c|cccccc}
		0 &  &   \\[3pt]
		\frac{1}{2} & \frac{1}{2}\varphi_{1,2}     \\[3pt]
		\frac{1}{2} & \frac{1}{2}\varphi_{1,3}-\varphi_{2,3} & \varphi_{2,3} \\[3pt]
		1 & \varphi_{1,4}-2\varphi_{2,4} & \varphi_{2,4} & \varphi_{2,4} \\[3pt]
		\frac{1}{2} & \frac{1}{2}\varphi_{1,5}-2a_{5,2}-a_{5,4} & a_{5,2} & a_{5,2} & \frac{1}{4}\varphi_{2,5}-a_{5,2}\\[3pt]
		\hline  & \varphi_1-3\varphi_2+4\varphi_3 & 0 & 0 & -\varphi_2+4\varphi_3 & 4\varphi_2-8\varphi_3
\end{array}\quad,\end{align}
with $a_{5,2}=\frac{1}{2}\varphi_{2,5}-\varphi_{3,4}+\frac{1}{4}\varphi_{2,4}-\frac{1}{2}\varphi_{3,5}$.
Although this five-stage method is fourth-order accurate for semilinear parabolic problems, it may be not a good candidate for solving the gradient system \eqref{problem: stabilized version}. 
Actually, it would not be stabilized to preserve the energy dissipation law \eqref{problem: energy dissipation law} unconditionally because the associated differentiation matrix $D^{(4,H)}(z)$ 
is not positive definite, see Figure \ref{fig: principal minors of EERK4 methods}(d), in which the curves of fourth and fifth leading principal minors of $\mathcal{S}(D^{(4,H)};z)$ are depicted.

Up to now, we are not able to find a fourth-order EERK method that preserves the energy dissipation law 
\eqref{problem: energy dissipation law} unconditionally. Nonetheless, this issue would be theoretically interesting and practically important in simulating the gradient system \eqref{problem: gradient flows}.

To end this article, we summarize our results in the following.
With a unified theoretical framework and a new indicator, namely average dissipation rate, for the energy dissipation properties of EERK methods, we examine some of popular methods and find: 
\begin{itemize}
	\item[(i)] Among second-order EERK methods, the average dissipation rate of the EERK2 method \eqref{scheme: EERK2 Butcher tableau} with $c_2=\tfrac{1}{2}$ is the closest to the continuous one so that it preserves the energy dissipation law \eqref{problem: energy dissipation law} best although the ETD2RK method \eqref{scheme: ETDRK2}, corresponding to the EERK2 method \eqref{scheme: EERK2 Butcher tableau} with $c_2=1$, seems the most popular for gradient flows, see \cite{CoxMatthews:2002JCP, DuJuLiQiao:2019SINUM, DuJuLiQiao:2021SIREV, JuLiQiaoZhang:2018MCOM,FuYang:2022JCP,LiuQuanWang:2023,ZhangLiuQianSong:2023}. If taking into the contractivity account, 
	the EERK2-w method \eqref{scheme: EERK2-weak Butcher tableau} with $c_2=\tfrac{1}{2}$ generates less time ``ahead" effect than the well-known ETD2RK method.
	\item[(ii)] Among third-order EERK methods, the popular ETD3RK and ETD2CF3 methods may destroy the energy dissipation law \eqref{problem: energy dissipation law}, especially for large time-step sizes. For the EERK3-1 \eqref{scheme: EERK3-HO1 Butcher} and EERK3-2 \eqref{scheme: EERK3-HO2 Butcher} methods, one can choose proper parameters (abscissas) to ensure the preserving of original dissipation law, while the EERK3-1 method \eqref{scheme: EERK3-HO1 Butcher} with $c_2=\tfrac{4}{9}$ produces the minimum time ``ahead" effect among the considered third-order EERK methods.
\end{itemize}
At the same time, our theory is far away from complete. There are many interesting issues that we have not yet addressed. Some of them are listed as follows:
\begin{itemize}
	\item[(a)] As mentioned, we are not able to find (or prove the non-existence of) a fourth-order EERK method that preserves the energy dissipation law \eqref{problem: energy dissipation law} unconditionally.	
		\item[(b)] It is noticed that the average dissipation rates $\mathcal{R}(z)$ of the mentioned EERK methods preserving the energy dissipation law \eqref{problem: energy dissipation law} are greater than 1 and unbounded, that is, $\mathcal{R}(z)\rightarrow+\infty$ as $z\rightarrow-\infty$. The method with a bounded average dissipation rate would be significantly preferred in the long-time adaptive simulation approaching the steady state. Is there such an EERK method or how do we construct it?
		\item[(c)] At least, is there a second-order EERK method that has a better  dissipation rate than the EERK2 method \eqref{scheme: EERK2 Butcher tableau} with $c_2=\tfrac{1}{2}$? Is there a third-order EERK method that has a better dissipation rate than the EERK3-1 method \eqref{scheme: EERK3-HO1 Butcher} with $c_2=\tfrac{4}{9}$?
\end{itemize}


	\section*{Acknowledgements}
The authors would like to thank Dr. Cao Wen for his sincere help.

\appendix

\section{Auxiliary functions for EERK3-1 methods \eqref{scheme: EERK3-HO1 Butcher}}

To examine the second and third leading principal minors 
of $\mathcal{S}(D^{(3,1)};c_2,z)$ with the differentiation 
matrix \eqref{def: EERK3-1 differentiation matrix}, we introduce 
two auxiliary functions $g_{31}$ and $g_{32}$ as follows,
{\scriptsize 	\begin{align}
	g_{31}(c_2,\sigma,z)&\,:=
	-9 c_2^2 z^2 e^{2 \sigma z}+18 c_2 z e^{(\sigma+\frac{4}{3}) z}
	+6 c_2 z e^{(\sigma+\frac{2}{3}) z} ((3 c_2-2) z-3)-9 e^{\frac{4 z}{3}}(c_2^2 z^2+1)\nonumber\\
	&\,	-6 e^{\frac{2 z}{3}} ((3 c_2-2) z-3)+(2 z+3) (2(3 c_2-1) z-3),\label{def: g31}\\
	g_{32}(c_2,\sigma,z)&\,:=
	-27 c_2^2 z(-z^2+2 z+3 e^{2 z}-2 e^z (z+3)+3)(e^{\frac{2 z}{3}}-e^{\sigma z})^2
	-4 (e^{2 z}-1) z(2 z-3 e^{\frac{2 z}{3}}+3)^2\nonumber\\
	&\,+6 c_2\kbra{-12 z^2 e^{(\sigma+\frac{4}{3}) z}+4 (2 z+3) z^2 e^{(\sigma+\frac{2}{3}) z}
		+2 (2 z^2+9 z+9) z e^{\sigma z+z}+18 z e^{(\sigma+\frac{8}{3}) z}
		-6 (z+3) z e^{(\sigma+\frac{5}{3}) z}}\nonumber\\
		&\,
		+6 c_2\kbra{-6 (2 z+3) z e^{(\sigma+2) z}+9 e^{\frac{2 z}{3}}(z^2-2 z-3)-6 e^z (2 z^2+9 z+9)
		-3 e^{\frac{8 z}{3}}(8 z^2+6 z+9)-6 e^{\frac{7 z}{3}} (z+3) z}\nonumber\\
		&\,
		+6 c_2\kbra{+2 e^{\frac{5 z}{3}}(2 z^3+9 z^2+18 z+27)+e^{2 z}(8 z^3+12 z^2+18 z+27)+3(-2 z^3+z^2+12 z+9)
		+18 e^{\frac{10 z}{3}} z}.\label{def: g32}
\end{align}}
For any constants $p_2>0$ and $p_1\ge0$, one has $\lim_{z\rightarrow-\infty}z^{p_1}e^{p_2z}\rightarrow0$. The dominant parts of $g_{31}$ and $g_{32}$ are simple although the expressions of them seem rather complex. Thus the computer-aided proof is always applied for simplicity of presentation.
We will prove the following results by applying the technique of comparison function developed in Propositions \ref{proposition: g21} and \ref{proposition: g22}. 
\begin{proposition}\label{proposition: g31}
 For the function $g_{31}$  in \eqref{def: g31}, $g_{31}(c_2,c_2,z)>0$ if $c_2\in[\tfrac4{9},1]$ and $z<0$.
\end{proposition}
\begin{proof}
	For the function $g_{31}$  in \eqref{def: g31}, we consider a comparison function $g^*_{31}(c_2,z)=g_{31}(c_2,\tfrac{4}{9},z)$ such that the difference 
		\begin{align*}
		g_{31}(c_2,c_2,z)-&\,g_{31}^*(c_2,z)=3 c_2(e^{\frac{4 z}{9}}-e^{c_2 z})e^{\frac{4 z}{9}}
		\braB{3 c_2  z^2+3 c_2 z^2 e^{(c_2-\frac{4}{9})z}-2e^{\frac{2 z}{9}} ((3 c_2-2) z^2-3z)-6z e^{\frac{8 z}{9}}}\\
		\ge&\,-3 c_2^2z(e^{\frac{4 z}{9}}-e^{c_2 z})e^{\frac{4 z}{9}}
		\braB{-3z \brat{1+e^{\frac{5z}{9}}-2e^{\frac{2 z}{9}}}+2c_2^{-1}e^{\frac{2 z}{9}}\brat{3e^{\frac{2 z}{3}}-3-2z}}\\
		\ge&\,-3 c_2^2z(e^{\frac{4 z}{9}}-e^{c_2 z})e^{\frac{4 z}{9}} r_{31}(z)\ge0\quad\text{for $c_2\in[\tfrac4{9},1]$ and $z<0$,}
	\end{align*}
	where the auxiliary function $r_{31}(z):=-3z \brat{1+e^{\frac{5z}{9}}-2e^{\frac{2 z}{9}}}+2e^{\frac{2 z}{9}}\brat{3e^{\frac{2 z}{3}}-3-2z}$ is decreasing and positive for $z<0$, cf. Figure \ref{fig: g31, g32 comparfuns} (a).
	Note that, $g^*_{31}(c_2,z)$ is a concave, quadratic polynomial with respect to $c_2$ because $$\partial_{c_2}^2g^*_{31}(c_2,z)=-9 e^{\frac{8 z}{9}} (e^{\frac{2 z}{9}}-1)^2 z^2<0.$$
	Through lengthy and simple calculations, it is not difficult to check that, cf. Figure \ref{fig: g31, g32 comparfuns} (b), 
	\begin{align*}
		g^*_{31}(1,z)>0\quad\text{and}
		\quad g^*_{31}(\tfrac{4}{9},z)>0\quad\text{for $z<0$.}
	\end{align*}
	They imply that $g_{31}^*(c_2,z)>0$ and then $g_{31}(c_2,c_2,z)>0$ for $c_2\in[\tfrac4{9},1]$ and $z<0$.
\end{proof}

\begin{figure}[htb!]
	\centering
	\subfigure[ $r_{31}(z)$]
	{\includegraphics[width=2.15in]{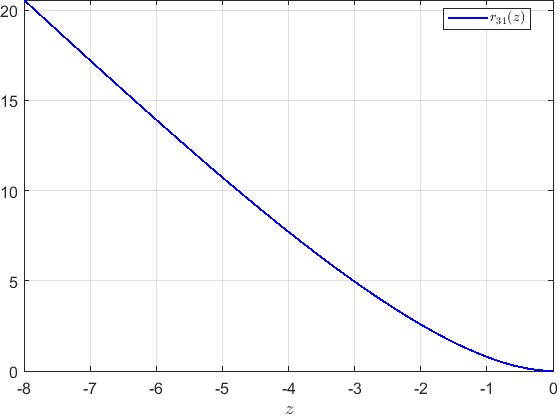}}\hspace{10mm}
	\subfigure[ $g^*_{31}(c_2,z)$]
	{\includegraphics[width=2.15in]{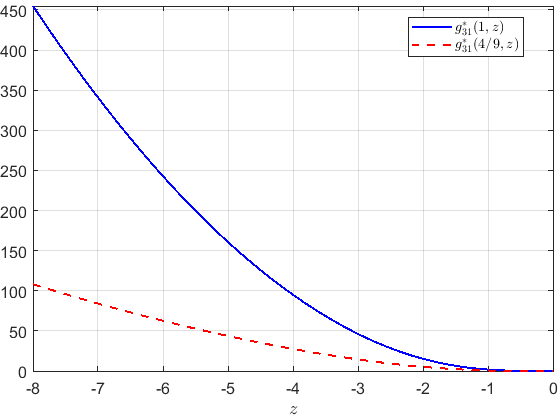}}\\
	\subfigure[ $r_{32,1}(z)$ and $r_{32,2}(z)$]
	{\includegraphics[width=2.15in]{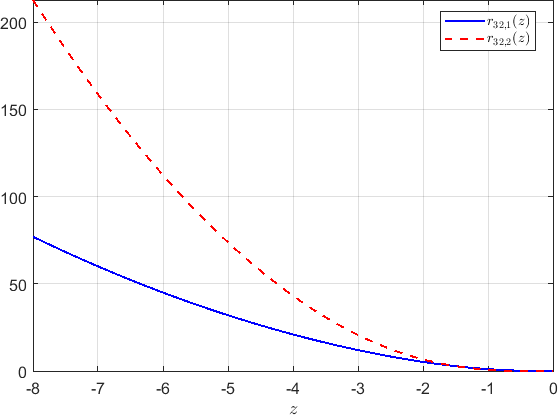}}\hspace{10mm}
	\subfigure[ $g^*_{32}(c_2,z)$]
	{\includegraphics[width=2.15in]{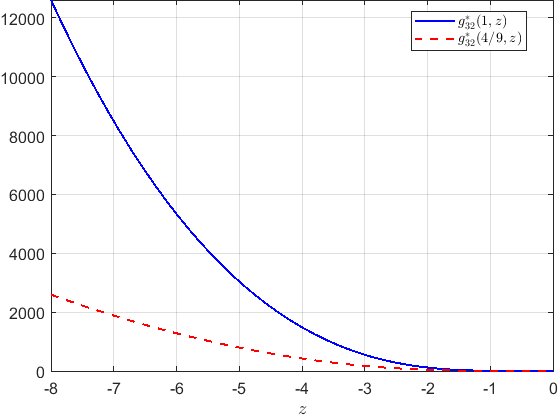}}
	\caption{Auxiliary functions $r_{31}$, $g^*_{31}$, $r_{32,1}$, $r_{32,2}$ and $g^*_{32}$.}
	\label{fig: g31, g32 comparfuns}
\end{figure}

\begin{proposition}\label{proposition: g32}
	For the function $g_{32}$  in \eqref{def: g32}, $g_{32}(c_2,c_2,z)>0$ if $c_2\in[\tfrac4{9},1]$ and $z<0$.
\end{proposition}
\begin{proof}
	For the function $g_{32}$  in \eqref{def: g32}, we consider a comparison function $g^*_{32}(c_2,z):=g_{32}(c_2,\tfrac{4}{9},z)$ such that the difference
		\begin{align*}
		g_{32}(c_2,c_2,z)-g_{32}^*(c_2,z)
		=&\,-3c_2^2 z(e^{\frac{4 z}{9}}-e^{c_2 z})e^{\frac{2z}{3}}
		\kbra{9  \brab{e^{(c_2-\frac2{3})z}+e^{-\frac{2z}{9}}-2}r_{32,1}(z)
		+4c_2^{-1}r_{32,2}(z)}\\
		\ge&\,-3c_2^2 z(e^{\frac{4 z}{9}}-e^{c_2 z})e^{\frac{2z}{3}}
		\kbra{9(e^{\frac{z}{3}}+e^{-\frac{2z}{9}}-2)r_{32,1}(z)
			+4r_{32,2}(z)},
	\end{align*}
	where the two auxiliary functions $r_{32,1}$ and  $r_{32,2}$ defined by
		\begin{align*}
		r_{32,1}(z):=&\,z^2-2 z-3 e^{2 z}+2 e^z (z+3)-3,\\
		r_{32,2}(z):=&\,4z^2+6z+e^{\frac{z}{3}}(2 z^2+9 z+9)-6 e^{\frac{2z}{3}} z+9 e^{2z}-3 e^{z} (z+3)-3 e^{\frac{4z}{3}} (2 z+3).
	\end{align*}
	 Since the functions $r_{32,1}$ and  $r_{32,2}$ are decreasing and positive for $z<0$, cf. Figure \ref{fig: g31, g32 comparfuns} (c), we see that $g_{32}(c_2,c_2,z)\ge g_{32}^*(c_2,z)$ for $c_2\in[\tfrac4{9},1]$ and $z<0$.

	Note that, $g^*_{32}(c_2,z)$ is a concave, quadratic polynomial with respect to $c_2$ due to $$\partial_{c_2}^2g^*_{32}
	=27z e^{\frac{8 z}{9}}(e^{\frac{2 z}{9}}-1)^2r_{32,1}(z)<0\quad\text{for $z<0$.}$$
	By simple but lengthy calculations, it is not difficult to check that, cf. Figure \ref{fig: g31, g32 comparfuns} (d), 
	\begin{align*}
		g^*_{32}(1,z)>0\quad\text{and}
		\quad g^*_{32}(\tfrac{4}{9},z)>0\quad\text{for $z<0$.}
	\end{align*}
	They imply that $g_{32}^*(c_2,z)>0$ and then $g_{32}(c_2,c_2,z)>0$ for $c_2\in[\tfrac4{9},1]$ and $z<0$.
\end{proof}

\section{Auxiliary functions for EERK3-2 methods \eqref{scheme: EERK3-HO2 Butcher}}

To examine the second and third leading principal minors 
of $\mathcal{S}(D^{(3,2)};1,\tfrac{1}{2},z)$, we introduce 
two auxiliary functions $g_{41}$ and $g_{42}$ as follows,
{\scriptsize\begin{align}
	g_{41}(z):=&\,
	5(3 z^2+2 z-5)+8 e^{\frac{3 z}{2}}(z^2-5 z-1)-2 e^z(8 z^2+z+3)\nonumber\\
	&\,+e^{2 z}(-16 z^2+32 z-1)-8 e^{z/2} (z-5)+8 e^{\frac{5 z}{2}} z.\label{def: g41}\\
	g_{42}(z):=&\,
		800+1241 z+334 z^2-71 z^3+8 e^{z/2} (z^2-97 z-80)-8 e^{\frac{5 z}{2}}(121 z^2+253 z+80)\nonumber\\
		&\,-2 e^{3 z}(149 z^2+19 z+80)+8 e^{\frac{3 z}{2}}(5 z^3+160 z^2+216 z+160)
		+2 e^z(40 z^3-331 z^2-981 z-880)\nonumber\\
		&\,+e^{2 z}(431 z^3+306 z^2+928 z+1120)+1072 e^{\frac{7 z}{2}} z-169 e^{4 z} z.\label{def: g42}
\end{align}}
For simplicity of presentation, we always use the computer-aided proof to prove the positivity of  the involved auxiliary functions.
\begin{figure}[htb!]
	\centering
	\subfigure[ $r_{41}(z)$]
	{\includegraphics[width=2.15in]{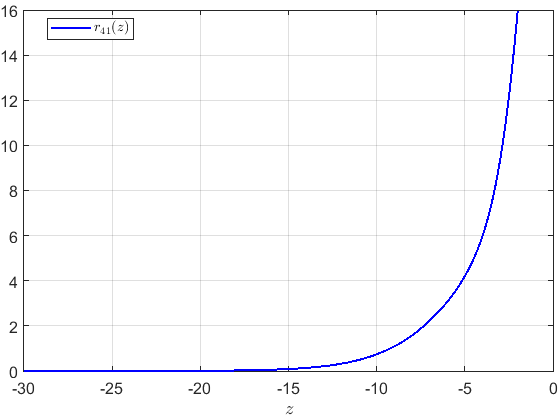}}\hspace{10mm}
	\subfigure[ $g_{41}(z)$]
	{\includegraphics[width=2.15in]{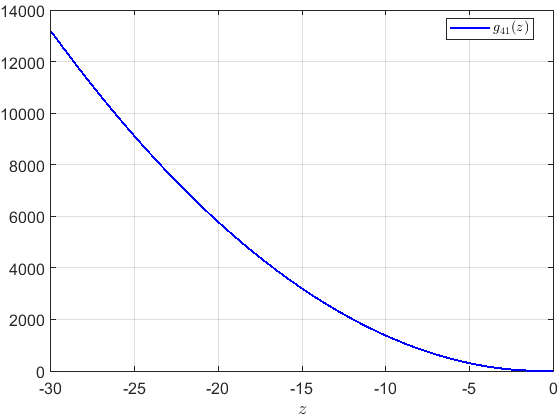}}\\
	\subfigure[ $r_{42}(z)$]
	{\includegraphics[width=2.15in]{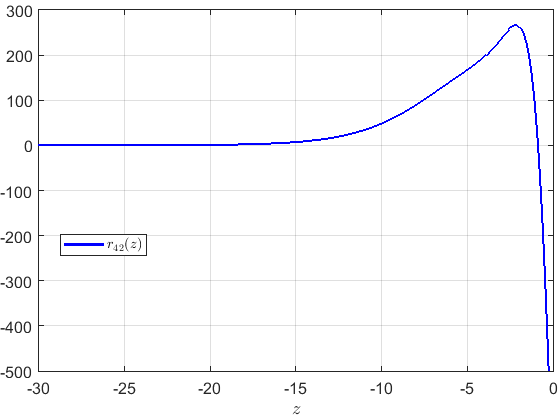}}\hspace{10mm}
	\subfigure[ $g_{42}(z)$]
	{\includegraphics[width=2.15in]{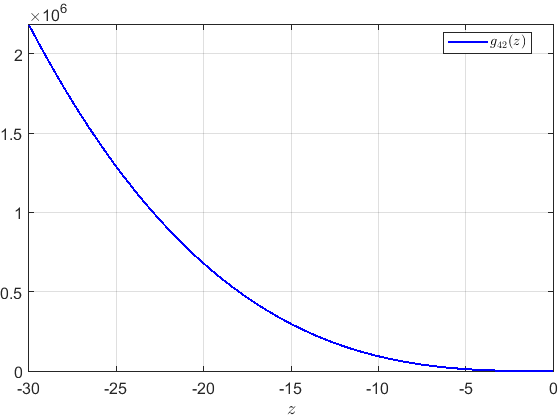}}\\
	\caption{Curves of the functions $g_{41}(z)$ and $g_{42}(z)$.}
	\label{fig: g41 g42 curves}
\end{figure}

\begin{proposition}\label{proposition: g41 g42}
	The functions $g_{41}$ and $g_{42}$ in \eqref{def: g41}-\eqref{def: g42} are positive for $z<0$.
\end{proposition}
\begin{proof}
	Note that, the quadratic polynomial part $\bar{g}_{41}(z):=5(3 z^2+2 z-5)$ of $g_{41}(z)$ is decreasing with respect to $z\in(-1/3,0)$ and $\lim_{z\rightarrow-\infty}\bar{g}_{41}(z)=+\infty$. The remaining part $r_{41}(z):=g_{41}(z)-\bar{g}_{41}(z)$ approaches zero when $\abs{z}$ is properly large such as $z\le z_0:=-30$, see Figure \ref{fig: g41 g42 curves}(a). Actually, $r_{41}(z_0)\approx-8.6\times10^{-5}$.  That is to say, $\bar{g}_{41}$ is dominant for $z\in(-\infty,z_0)$. As seen in Figure \ref{fig: g41 g42 curves}(b), $g_{41}$ is decreasing and positive inside the finite interval $(z_0,0)$. They lead to $g_{41}(z)>0$ for $z<0$.
	
	Similarly, the cubic polynomial part $\bar{g}_{42}(z):=800+1241 z+334 z^2-71 z^3$ of $g_{42}(z)$ is decreasing for $z\in(-3/2,0)$ and $\lim_{z\rightarrow-\infty}\bar{g}_{42}(z)=+\infty$. The remaining part $r_{42}(z):=g_{42}(z)-\bar{g}_{42}(z)$ approaches zero when $\abs{z}$ is properly large such as $z\le z_0:=-30$, see Figure \ref{fig: g41 g42 curves}(c). Actually, $r_{42}(z_0)\approx-9.1\times10^{-3}$.  That is to say, $\bar{g}_{42}$ is dominant for $z\in(-\infty,z_0)$. As seen in Figure \ref{fig: g41 g42 curves}(d), $g_{42}$ is decreasing and positive inside $(z_0,0)$. They imply that $g_{42}(z)>0$ for $z<0$.
\end{proof}

To examine the second and third leading principal minors 
of $\mathcal{S}(D^{(3,2)};\frac{3}{4},\frac{3}{5},z)$, we define
two auxiliary functions $g_{51}$ and $g_{52}$ as follows,
{\scriptsize 	
	\begin{align}
	g_{51}(z)&\,:=\frac1{6400}\left[8(567 z^2-1353 z-3362)
	-625 e^{\frac{6 z}{5}}(9 z^2+16)-50 e^{\frac{27 z}{20}}(99 z^2+492 z+256)\right.\nonumber\\
	&\,\left.-e^{\frac{3 z}{2}}(5625 z^2+4096)+9600 e^{\frac{21 z}{10}} z+15000 e^{\frac{39 z}{20}} z
	+128 e^{\frac{3 z}{4}} (33 z+164)+200 e^{\frac{3 z}{5}} (33 z+164)\right],\label{def: g51}
\\
		g_{52}(z)&\,:=\frac1{8\times10^6}\left[8(-96681 z^3+600979 z^2+1843541 z+1050625)
		+1875 e^{\frac{6 z}{5}}(533 z^2-1148 z-385) z-2332800 e^{\frac{21 z}{10}} z^2\right.\nonumber\\
		&\,+150 e^{\frac{27 z}{20}}(7319 z^2+46552 z+46361) z+3 e^{\frac{3 z}{2}}(453125 z^2+62500 z-58849) z-1000 e^{\frac{13 z}{5}}(2754 z^2+1107 z+5125)\nonumber\\
		&\,-400 e^z (6018 z^2+35569 z+42025)-128 e^{\frac{3 z}{4}}(7319 z^2+46552 z+25625)-200 e^{\frac{3 z}{5}}(7319 z^2+46552 z+25625)\nonumber\\
		&\,-8 e^{\frac{11 z}{4}}(529821 z^2+54243 z+410000)+200 e^{\frac{8 z}{5}}(2106 z^3+37287 z^2+52087 z+51250)-3645000 e^{\frac{39 z}{20}} z^2\nonumber\\
		&\,+200 e^{\frac{7 z}{4}}(6642 z^3+24111 z^2+31963 z+32800)+8 e^{2 z}(137781 z^3-328779 z^2-65091 z+1050625)-7203 e^{\frac{7 z}{2}} z\nonumber\\
		&\,\left.-46875 e^{\frac{16 z}{5}} z+7650750 e^{\frac{67 z}{20}} z-7350 e^{\frac{5 z}{2}} (2 z-25) z+18750 e^{\frac{11 z}{5}} (14 z+41) z-300 e^{\frac{47 z}{20}} (6482 z+48683) z\right].\label{def: g52}
\end{align}}

\begin{figure}[htb!]
	\centering
	\subfigure[ $r_{51}(z)$]
	{\includegraphics[width=2.15in]{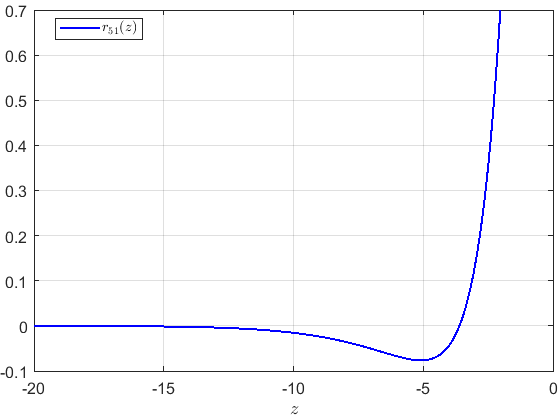}}\hspace{10mm}
	\subfigure[ $g_{51}(z)$]
	{\includegraphics[width=2.15in]{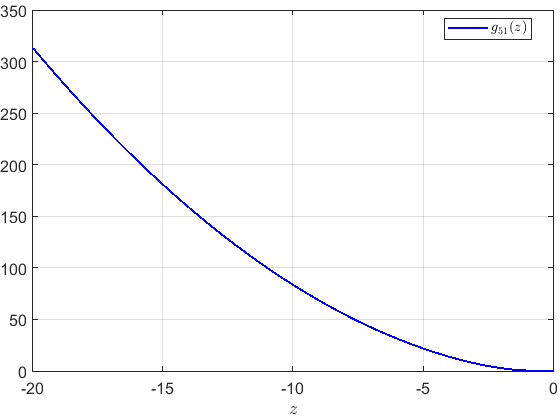}}\\
	\subfigure[ $r_{52}(z)$]
	{\includegraphics[width=2.15in]{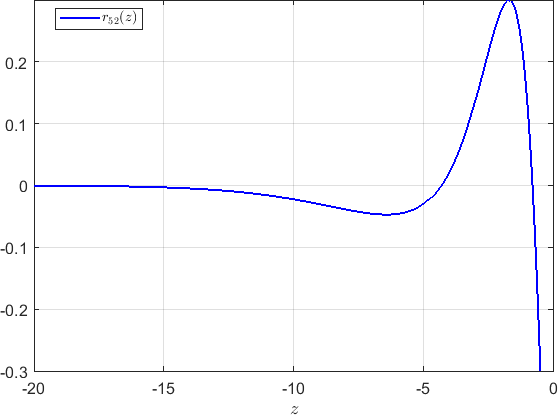}}\hspace{10mm}
	\subfigure[ $g_{52}(z)$]
	{\includegraphics[width=2.15in]{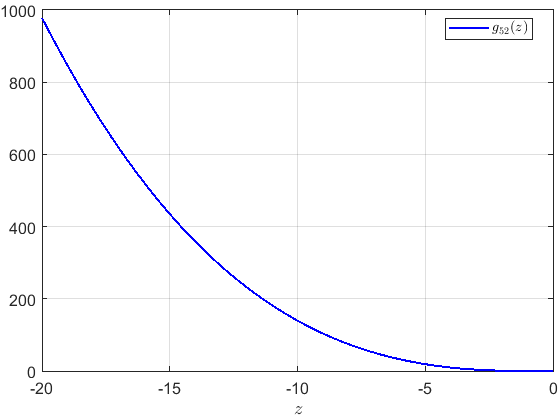}}\\
	\caption{Curves of the functions $g_{51}(z)$ and $g_{52}(z)$.}
	\label{fig: g51 g52 curves}
\end{figure}

\begin{proposition}\label{proposition: g51 g52}
	The functions $g_{51}$ and $g_{52}$ in \eqref{def: g51}-\eqref{def: g52} are positive for $z<0$.
\end{proposition}
\begin{proof}
	The quadratic polynomial part $\bar{g}_{51}(z):=\frac{1}{800}(567 z^2-1353 z-3362)$ of $g_{51}(z)$ is decreasing with respect to $z\in(-\infty,0)$ and $\lim_{z\rightarrow-\infty}\bar{g}_{51}(z)=+\infty$. The remaining part $r_{51}(z):=g_{51}(z)-\bar{g}_{51}(z)$ approaches zero when $\abs{z}$ is properly large such as $z\le z_0:=-20$, see Figure \ref{fig: g51 g52 curves}(a). Actually, $r_{51}(z_0)\approx-9.8\times10^{-5}$.  That is to say, $\bar{g}_{51}$ is dominant for $z\in(-\infty,z_0)$. As seen in Figure \ref{fig: g51 g52 curves}(b), $g_{51}$ is decreasing and positive inside $(z_0,0)$. It is easy to conclude that $g_{51}(z)>0$ for $z<0$.
	
	Similarly, the cubic polynomial part $\bar{g}_{52}(z):=\frac1{10^6}(-96681 z^3+600979 z^2+1843541 z+1050625)$ of $g_{52}(z)$ is decreasing for  $z\in(-3/2,0)$ and $\lim_{z\rightarrow-\infty}\bar{g}_{52}(z)=+\infty$. The remaining part $r_{52}(z):=g_{52}(z)-\bar{g}_{52}(z)$ approaches zero when $\abs{z}$ is properly large such as $z\le z_0:=-20$, see Figure \ref{fig: g51 g52 curves}(c). Actually, $r_{52}(z_0)\approx-3.2\times10^{-4}$.  That is to say, $\bar{g}_{52}$ is dominant for $z\in(-\infty,z_0)$. As seen in Figure \ref{fig: g51 g52 curves}(d), $g_{52}$ is decreasing and positive inside $(z_0,0)$. They imply that $g_{52}(z)>0$ for $z<0$.
\end{proof}

To examine the second and third leading principal minors 
of $\mathcal{S}(D^{(3,2)};\tfrac{1}{2},7/10,z)$, we define
two auxiliary functions $g_{61}$ and $g_{62}$ as follows,
{\scriptsize 	
	\begin{align}
		g_{61}(z)&\,:=\frac{1}{16}\left[-625 e^{\frac{7 z}{5}}(z^2+4)-50 e^{\frac{6 z}{5}}(17 z^2+64 z+28)+16 (21 z^2-136 z-256)\right.\nonumber\\
		&\,\left.-e^z(625 z^2+196)+700 e^{\frac{17 z}{10}} z+2500 e^{\frac{19 z}{10}} z+28 e^{z/2} (17 z+64)+100 e^{\frac{7 z}{10}} (17 z+64)\right],\label{def: g61}\\
		g_{62}(z)&\,:=\frac{1}{2000000}\left[
		2560000+4806561 z+2007966 z^2-25151 z^3-1102500 e^{\frac{19 z}{10}} z^2+2500 e^{\frac{7 z}{5}} \left(80 z^2-176 z+185\right) z\right.\nonumber\\
		&\,+50 e^{\frac{6 z}{5}} \left(8597 z^2+38299 z+32348\right) z-250 e^{\frac{27 z}{10}} \left(1911 z^2-2688 z+8000\right)-28 e^{z/2} \left(8597 z^2+38299 z+20000\right)\nonumber\\
		&\,-100 e^{\frac{7 z}{10}} \left(8597 z^2+38299 z+20000\right)-14 e^{\frac{5 z}{2}} \left(76881 z^2+63552 z+40000\right)\nonumber\\
		&\,+350 e^{\frac{3 z}{2}} \left(1029 z^3+4443 z^2+5606 z+3200\right)+50 e^{\frac{17 z}{10}} \left(6027 z^3+33159 z^2+63158 z+80000\right)\nonumber\\
		&\,+e^{2 z} \left(148176 z^3-893416 z^2+318039 z+2560000\right)+e^z \left(265625 z^3-1032550 z^2-5015039 z-5120000\right)\nonumber\\
		&\,\left.+2165500 e^{\frac{16 z}{5}} z-62500 e^{\frac{17 z}{5}} z-109561 e^{3 z} z-12500 e^{\frac{12 z}{5}} (11 z+32) z-50 e^{\frac{11 z}{5}} (17509 z+75658) z\right].\label{def: g62}
\end{align}}

\begin{figure}[htb!]
	\centering
	\subfigure[ $r_{61}(z)$]
	{\includegraphics[width=2.15in]{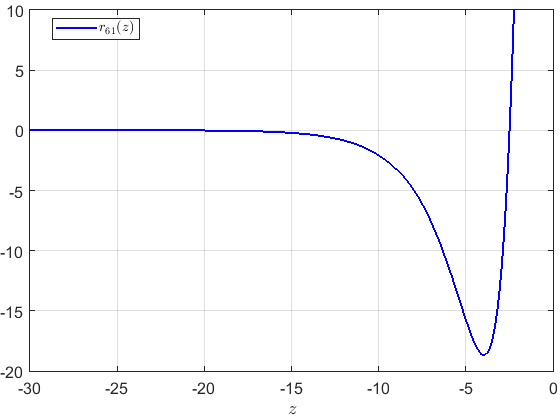}}\hspace{10mm}
	\subfigure[ $g_{61}(z)$]
	{\includegraphics[width=2.15in]{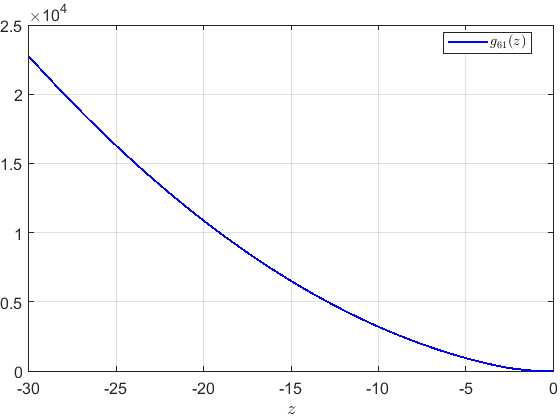}}\\
	\subfigure[ $r_{62}(z)$]
	{\includegraphics[width=2.15in]{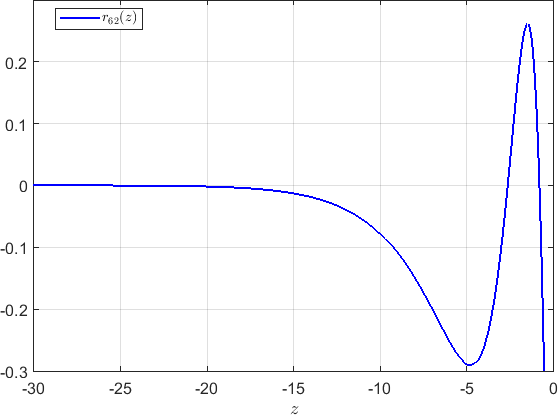}}\hspace{10mm}
	\subfigure[ $g_{62}(z)$]
	{\includegraphics[width=2.15in]{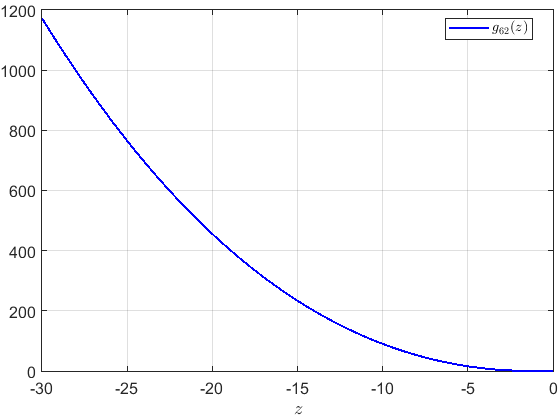}}\\
	\caption{Curves of the functions $g_{61}(z)$ and $g_{62}(z)$.}
	\label{fig: g61 g62 curves}
\end{figure}

\begin{proposition}\label{proposition: g61 g62}
	The functions $g_{61}$ and $g_{62}$ in \eqref{def: g61}-\eqref{def: g62} are positive for $z<0$.
\end{proposition}
\begin{proof}
	The quadratic polynomial part $\bar{g}_{61}(z):=21 z^2-136 z-256$ of $g_{61}(z)$ is decreasing with respect to $z\in(-\infty,0)$ and $\lim_{z\rightarrow-\infty}\bar{g}_{61}(z)=+\infty$. The remaining part $r_{61}(z):=g_{61}(z)-\bar{g}_{61}(z)$ approaches zero when $\abs{z}$ is properly large such as $z\le z_0:=-30$, see Figure \ref{fig: g61 g62 curves}(a). Actually, $r_{61}(z_0)\approx-2.4\times10^{-4}$.  That is to say, $\bar{g}_{61}$ is dominant for $z\in(-\infty,z_0)$. As seen in Figure \ref{fig: g61 g62 curves}(b), $g_{61}$ is decreasing and positive inside $(z_0,0)$. It is easy to conclude that $g_{61}(z)>0$ for $z<0$.
	
	Similarly, the cubic polynomial part $\bar{g}_{62}(z):=\frac{1}{2000000}\brat{
	2560000+4806561 z+2007966 z^2-25151 z^3}$ of $g_{62}(z)$ is decreasing for  $z\in(-3/2,0)$ and $\lim_{z\rightarrow-\infty}\bar{g}_{62}(z)=+\infty$. The remaining part $r_{62}(z):=g_{62}(z)-\bar{g}_{62}(z)$ approaches zero when $\abs{z}$ is properly large such as $z\le z_0:=-30$, see Figure \ref{fig: g61 g62 curves}(c). Actually, $r_{62}(z_0)\approx-2.9\times10^{-5}$.  That is to say, $\bar{g}_{62}$ is dominant for $z\in(-\infty,z_0)$. As seen in Figure \ref{fig: g61 g62 curves}(d), $g_{62}$ is decreasing and positive inside $(z_0,0)$. They imply that $g_{62}(z)>0$ for $z<0$.
\end{proof}

\end{document}